\documentclass[12pt,a4paper]{amsart}

\usepackage{amsmath, amsfonts, xifthen, latexsym, amssymb, amsthm, amscd}

\usepackage[utf8]{inputenc}
\usepackage{graphicx}
\usepackage{url}

\usepackage{hyperref}
\hypersetup{colorlinks=true,citecolor=blue,filecolor=blue,linkcolor=blue,urlcolor=blue}
\usepackage[margin=1.4in]{geometry}
\usepackage{fancyhdr}
\newtheorem{theorem}{Theorem}
\newtheorem{lemma}{Lemma}[section]

\newtheorem{corollary}{Corollary}[section]

\theoremstyle{definition}
\newtheorem{definition}{Definition}[section]

\newtheorem{remark}{Remark}

\newtheorem{example}{Example}
\newcommand{\uN}{{\mathbf{N}}}
\newcommand{\uM}{{\mathbf{M}}}
\newcommand{\R}{{\mathbb{R}}}
\newcommand{\N}{{\mathbb{N}}}
\newcommand{\RR}{{\R^\infty}}
\newcommand{\uP}{{\mathbf{P}}}
\newcommand{\hpi}{{\hat{\pi}}}
\newcommand{\hmu}{{\hat{\mu}}}
\newcommand{\halpha}{{\hat\alpha}}
\newcommand{\e}{\varepsilon}
\newcommand{\uZ}{{\mathbf{Z}}}
\newcommand{\uB}{{\mathbf{B}}}
\newcommand{\bB}{{\hat{\uB}}}
\newcommand{\ul}{{\mathbf{l}}}
\newcommand{\bl}{{\hat{\ul}}}
\newcommand{\uf}{{\mathbf{f}}}

\newcommand{\uy}{{\mathbf{y}}}
\newcommand{\ud}{{\mathbf{d}}}
\newcommand{\bd}{{\hat{\ud}}}
\newcommand{\ua}{{\mathbf{a}}}
\newcommand{\ux}{{\mathbf{x}}}

\newcommand{\ba}{{\hat{\ua}}}
\newcommand{\ub}{{\mathbf{b}}}
\newcommand{\bb}{{\hat{\ub}}}
\newcommand{\up}{{\mathbf{p}}}
\newcommand{\bp}{{\hat{\up}}}
\newcommand{\bx}{{\hat{\ux}}}
\newcommand{\by}{{\hat{\uy}}}
\newcommand{\uq}{{\mathbf{q}}}
\newcommand{\bq}{{\hat{\uq}}}

\newcommand{\uz}{{\mathbf{z}}}
\newcommand{\bz}{{\hat{\uz}}}
\newcommand{\uC}{{\mathbf{C}}}
\newcommand{\uF}{{\mathbf{F}}}
\newcommand{\bF}{{\hat{\uF}}}
\newcommand{\uH}{{\mathbf H}}
\newcommand{\uX}{{\mathbf{X}}}

\newcommand{\bX}{{\hat{\uX}}}

\newcommand{\uA}{{\mathbf{A}}}
\newcommand{\uS}{{\mathbf{S}}}
\newcommand{\bS}{{\hat{\uS}}}
\newcommand{\bA}{{\hat{\uA}}}

\newcommand{\cA}{{\mathcal{A}}}

\newcommand{\cP}{{\mathcal{P}}}
\newcommand{\cM}{{\mathcal{M}}}

\newcommand{\diam}{\mathrm{diam}}
\newcommand{\core}{{\mathrm{Core}}}

\newcommand{\prob}{\mbox{Prob}}
\newcommand{\bbf}{{\hat{\uf}}}

\newcommand{\limo}{{\lim_\omega}}
\newcommand{\prodo}{{\prod_\omega}}
\newcommand{\hprodo}{{{\hat{\prod}}_\omega}}
\newcommand{\uT}{{\mathbf{T}}}
\newcommand{\bT}{{\hat{\uT}}}

\newcommand{\bC}{{\hat{\uC}}}
\newcommand{\uD}{{\mathbf{D}}}
\newcommand{\bD}{{\hat{\uD}}}
\newcommand{\bH}{{\hat{\uH}}}

\begin{document}
\title[Limits of finite trees]{Convergence and limits of finite trees}
\author{Gábor Elek$^{a,b}$ and Gábor Tardos$^{b,c}$\\
\\
\small $^a$ Department of Mathematics and Statistics, Fylde College,\\Lancaster
University, Lancaster, United Kingdom\\
\small $^b$ Alfr\'ed R\'enyi Institute of Mathematics, Budapest, Hungary\\
\small $^c$ Moscow Institute of Physics and Technology, Dolgoprudny, Russia\\
\\
\small \texttt{g.elek@lancaster.ac.uk, tardos@renyi.hu}
}  

\thanks{The first author was partially supported
by the ERC Consolidator Grant ``Asymptotic invariants of discrete groups,
sparse graphs and locally symmetric spaces'' No. 648017 and by the ERC Synergy
Grant ``Dynasnet'' No. 810115.\\
\indent
The second author was partially supported by the ERC Synergy Grant
``Dynasnet'' No. 810115, the ERC advanced grant ``GeoSpace'' No. 882971, the National Research, Development and Innovation Office — NKFIH
projects K-116769, K-132696, KKP-133864, SNN-117879, SSN-135643 and by the grant of Russian Government N 075-15-2019-1926.}

\begin{abstract}
Motivated by the work of Lov\'asz and Szegedy on the convergence and limits of dense graph sequences \cite{LSZ}, we investigate the convergence and limits of finite trees with respect to sampling in normalized distance. We introduce dendrons (a notion based
on separable real trees) and show that the sampling limits of
finite trees are exactly the dendrons. We also prove that the limit dendron is unique.
\end{abstract}

\begingroup
\let\MakeUppercase\relax
\maketitle
\endgroup


\section{Introduction}
The main motivation of our paper is the (sampling) limit theory of dense graphs
introduced by Lov\'asz and Szegedy \cite{LSZ}. Let us recall very briefly the 
most important definitions. Let $G$ be a finite simple graph and $r\ge1$ be an
integer. Let us pick $r$ distinct vertices $v_1,v_2,\dots v_r$ of $G$ uniformly at random
and consider the graph $H$ induced by the chosen vertices. (Here we assume that the graph has at least $r$ vertices. Alternatively, one can use sampling with repetition and define $H$ as the $r$ vertex graph obtained by possible duplication of the vertices in the induced subgraph.) Then, $H$ will be
isomorphic to one of the $2^{r \choose 2}$ graphs on $r$ labeled vertices. Thus, the random choice of the vertices $v_i$
defines a probability distribution $p^G_r$ on the finite set $A_r$ of these labeled graphs. We say that
the sequence of finite graphs $(G_n)_{n\in\N}$ is {\bf convergent} if
$\lim_{n\to\infty} p^{G_n}_r(K)$ exists for all $r\in\N$ and $K\in A_r$. (Throughout this paper $\N$ stands for the set of positive integers.) Lov\'asz and Szegedy constructed a universal limit object for such
convergent graph sequences, the {\bf graphons}. A graphon is a measurable symmetric
function $W:[0,1]^2\to [0,1]$. For $K\in A_r$ (with the vertices of $K$ denoted by the integers $1$ through $r$),
$p^W_r(K)$ is defined as
$$p^W_r(K)=\int_0^1\int_0^1\dots\int_0^1  \prod_{i<j,(i,j)\in E(K)}
W(x_i,x_j)\prod_{i<j,(i,j)\notin E(K)}
(1-W(x_i,x_j))dx_1 dx_2\dots dx_r\,.$$
The graphon $W$ is the limit of the sequence $(G_n)_{n\in\N}$ if
for all $r\geq 1$ and $K\in A_r$, $\lim_{n\to\infty} p^{G_n}_r(K)= p^W_r(K)$.
It has been proved in \cite{LSZ} that for any convergent sequence
 $(G_n)_{n\in\N}$ there exists a limit graphon and all graphons
are limits of convergent sequences of finite graphs. The uniqueness problem
was considered in \cite{Borgs}.
The goal of our paper is to introduce and study the sampling limit theory
of finite trees. This aim seems to be rather contradictory since the
sampling limit theory of Lov\'asz and Szegedy gives non-trivial limit objects
only if the graphs are dense, and the trees are very sparse graphs.
We solve this problem by regarding trees as dense objects using their natural metric structure.

We identify each finite tree with its vertex
set. We call a finite tree {\bf nontrivial} if it has at least two vertices. Let $T$ be a nontrivial finite tree. In order to have more sequences of finite trees that are convergent we need to normalize the shortest path metric $\overline{d}_T$ on $T$. Here we choose to normalize it by its diameter $\diam_{\overline{d}_T}(T)$. We consider the following metric space structure on $T$.
$$d_T(x,y)= \frac{\overline{d}_T(x,y)}{\diam_{\overline{d}_T}(T)}.$$
Throughout the paper we use this normalization making every non-trivial finite tree into a metric space of diameter exactly 1. This is for convenience only; see Remark~\ref{remarklong} on the possibility for other normalizations. We consider the
uniform probability measure $\mu_T$ on $T$ as well. So, we obtain a 
{\bf metric measure space} structure on our tree. The notion of sampling convergence
for  metric measure spaces has been introduced by Gromov in the famous Section
$3\frac{1}{2}$ of his treatise ``Metric structures for Riemannian and
Non-Riemannian Spaces'' \cite{Gromov}. Let us recall the formal definitions.

\begin{definition}\label{sampling}
For $r\geq 1$, let $M_r$ be the space of real matrices $(d_{ij})_{1\le i,j\le r}$. We identify $M_r$ with the product space  $\prod_{1\le i,j\le r}\R$.

For a set $X$ and a function $d:X^2\to\R$ we define the map $\rho_r=\rho_r^{X,d}:X^r\to M_r$ with $\rho_r(x_1,x_2,\dots,x_r)=(d_{ij})_{1\le i,j\le r}$, where $d_{ij}=d(x_i,x_j)$ for $i\ne j$ and $d_{ii}=0$. Treating $d_{ii}$ separately is needed because we will later use $\rho_r$ (and the sampling measure $\tau_r$ below) not only for distance functions $d$ but also for functions not satisfying $d(x,x)=0$ for all $x\in X$.

A {\bf metric measure space} is a triple $(X,d,\mu)$, where $(X,d)$ is a complete separable metric space and $\mu$ is a probability Borel measure on $X$.

For a metric measure space, or more generally for a triple $(X,d,\mu)$, where $\mu$ is a probability measure on $X$ and $d:X^2\to\R$ is $\mu^2$-measurable we define the sampling measure $\tau_r(X,d,\mu)\in\prob(M_r)$ as the push-forward of the product measure $\mu^r$ along $\rho_r$, that is, for a Borel set $H$ in $M_r$ we set $\tau_r(X,d,\mu)(H)=\mu^r(\rho_r^{-1}(H))$.

For a nontrivial finite tree $T$ and $r\ge1$ we write $\tau_r(T)$ as a shorthand for $\tau_r(T,d_T,\mu_T)$, where $d_T$ is the normalized distance and $\mu_T$ is the uniform measure. Our normalization makes the diameter $1$, so the image of $\rho_r^{T,d_T}$ is contained in $M_r^1=\prod_{1\le i,j\le r}[0,1]$ and thus $\tau_r(T)$ is concentrated on $M_r^1$.
\end{definition}

\begin{definition}\label{converge}
The sequence of finite trees $(T_n)_{n\in\N}$ is
convergent in metric sampling if for any $r\geq 1$ the sequence of
measures $(\tau_r(T_n))_{n\in\N}$ converges weakly in the space $\prob(M_r)$ of the Borel probability measures on $M_r$.
\end{definition}

Since $\tau_r(T)$ is concentrated on the compact space $M_r^1$ for any nontrivial finite tree $T$, one can pick a convergent subsequence from any sequence
of finite trees. The main goal of our paper is to identify the limit objects
of such convergent sequences of finite trees. First, we need to recall the
notion of a real tree, our key topological notion (see
\cite{Bestvina} for a survey).

\begin{definition}\label{realtree}
We say that the non-empty complete metric space $(T,d)$ is a {\bf real tree} or {\bf$\R$-tree} if for any pair of distinct points $p,q\in T$ one has an isometric embedding of an interval $\alpha:[a,b]\to X$ such that $\alpha(a)=p$, $\alpha(b)=q$ and $\alpha(c)$ separates $p$ from $q$ for any $a<c<b$, that is, $p$ and $q$ are in distinct connected components of $T\setminus\{\alpha(c)\}$.

We write $[x,y]$ to denote $\alpha([a,b])$, which is uniquely determined in a real tree. We let $[x,x]=\{x\}$ and call $[x,y]$ a {\bf segment} in the real tree $(T,d)$. If $x\ne y$ we call $[x,y]$ a {\bf proper segment}. All the points of $[x,y]$ other than its end points $x$ and $y$ are called the {\bf intermediate points} of the segment $[x,y]$. An {\bf intermediate point} of a real tree is an intermediate point of a proper segment in the tree.

For $p\in T$ we call the connected components of $T\setminus\{p\}$ the {\bf $p$-branches}. Note that $p$ is an intermediate point if and only if there are at least two $p$-branches. A {\bf branch} of $T$ is a $p$-branch for some $p\in T$.
\end{definition}

A natural limit object would be a {\bf measured real tree}, that is a separable real tree equipped with a probability measure making it a metric measure space. In fact, there are several known metrics for metric measure spaces, like the Gromov-Prohorov metric \cite{GPW} and Gromov's
$\underline{\square}_1$ metric \cite{Gromov}. Both of these define Gromov's
weak topology (given by the sampling) \cite{Lohr}, see also Section 10 of \cite{Jans}. This means that a sequence of measured real trees tend to another metric measure space (necessarily a measured real tree) in one of these metrics if and only if they tend to it in metric sampling. Nevertheless, these metrics are substantively different from the metric obtained from sampling that we consider in this paper. In particular,
Gromov-Prohorov metric and the $\underline{\square}_1$ metric is complete as opposed to the metric obtained from metric samplings, which is pre-compact. 

We will introduce \emph{dendrons} as limit objects for sequences of finite trees that are convergent in metric sampling. If a sequence of finite trees (considered as metric measure spaces after normalization) tend to another metric measure space in Gromov's weak topology, then the finite trees form a convergent sequence also in metric sampling and the limit dendron can be easily obtained from the limit in the Gromov's weak topology, see Example~1 below. But there are many sequences of finite trees that are convergent in metric sampling, but the corresponding measured real trees are not convergent in the Gromov's weak topology, see examples~2--6 below.
Hence, our dendrons can be considered to be the compactification of the space
of measured real trees of maximum diameter $1$.

It is natural to construct the limit objects as some generalization of measured real trees. One way to do it is to drop the requirement of a measured real tree to be separable (and at the same time allow some Borel sets not to be measurable). See the definition of quasi measured real trees in the next section and Remarks~\ref{meareal} and \ref{strong}. Here we take another route.

\begin{definition}\label{dendron}
A {\bf long dendron} $D=(T,d,\nu)$ is a real tree $(T,d)$ together with a probability Borel measure $\nu$ on $A_D=T\times[0,\infty)$ satisfying $\nu(B\times[0,\infty))>0$ for all branches $B$ of $T$. We define $d_D:A_D^2\to\R$ by $d_D((u,a),(v,b))=d(u,v)+a+b$. We say that $D$ is a {\bf dendron} if $d_D$ is almost surely bounded by $1$, that is, if $\nu^2(\{(x,y)\in A_D^2\mid d_D(x,y)>1\})=0$.

For a long dendron $D=(T,d,\nu)$ and $r\ge1$, we write $\tau_r(D)$ as the shorthand for the sampling measure $\tau_r(A_D,d_D,\nu)$.

We will use the notation $A_D$ and $d_D$ for (long) dendrons $D$ in the above sense. We say that two long dendrons $D=(T,d,\nu)$ and $D'=(T',d',\nu')$ are {\bf isomorphic} if there exists an isometry $f$ from $(T,d)$ to $(T',d')$ such that $f':A_D\to A_{D'}$ defined as $f'(p,a)=(f(p),a)$ is measure preserving.
\end{definition}

\begin{remark} 
In the definition of a (long) dendron we did not explicitly require the real tree $(T,d)$ to be separable, but this is implicit in the definition. We will state this explicitly in Corollary~\ref{corsep}.

Note that $d_D$ defined above for a (long) dendron $D=(T,d,\nu)$ is not a distance on $A_D$ as $d_D(x,x)=2a>0$ for all $x=(u,a)\in A_D$ with $a>0$.

Long dendrons can formally be considered marked metric measure spaces as introduced by Depperschmidt, Greven and Pfaffelhuber, see \cite{DGP}, where the infinite interval $[0,\infty)$ is the space of possible marks. This connection is superficial though, as the topology is rather different. We will see that finite trees correspond to dendrons where the mark is $0$ with probability one, but their limit points include all dendrons, among them those where the mark is separated from $0$ with probability one. See Examples~2, 3 and 6 below.
\end{remark}

Our main results are the following three theorems that are analogues of the main
results of Lov\'asz and Szegedy in \cite{LSZ} and Borgs, Chayes and Lov\'asz in \cite{Borgs}.

\begin{theorem}\label{tetel1}
For any convergent sequence of finite trees $(T_n)_{n\in\N}$ there
exists a dendron $D$ (the sampling limit of $(T_n)_{n\in\N}$) such that the sampling measures $\tau_r(T_n)$ weakly converge to $\tau_r(D)$ for all positive integers $r$.
\end{theorem}

\begin{theorem}
\label{tetel2}
Any dendron is the sampling limit of a convergent sequence of finite trees.
\end{theorem}

\begin{theorem}\label{unique}
The sampling limit is unique up to isomorphism. More generally, if two long dendrons $D$ and $D'$ satisfy $\tau_r(D)=\tau_r(D')$ for all $r$,  then $D$ and $D'$ are isomorphic.
\end{theorem}

\begin{remark}\label{remarklong}
These theorems establish dendrons as the sampling limits of finite trees with respect to distance normalized by the diameter. Long dendrons can be considered an extension suitable for other normalizations. It is possible to consider finite trees with explicitly given normalization factors. In this case they correspond to metric measure spaces with an arbitrary diameter. We formulate Theorem~\ref{unique} for long dendrons to be more general and capture those limits with an unbounded diameter too. The single main difference between this more general limit theory of finite trees and the theory discussed in this paper is that with no bound on the diameter we lose compactness: we will not be able to find convergent subsequences of any sequence of arbitrarily normalized finite trees. Janson in a recent paper, \cite{Jans} found other normalization factors more suitable in various (often probabilistic) scenarios.

For the purposes of this paper we encourage the reader to concentrate on the case of dendrons. For this case it is instructive to note that one could define the domain $A_D$ of a dendron $D=(T,d,\nu)$ as $A_D=T\times[0,1/2]$ instead of the definition $A_D=T\times[0,\infty)$ above. Indeed, for the complementary set $B=T\times(1/2,\infty)$ we have $d_D(x,y)>1$ for all $x,y\in B$, hence if $D$ is a dendron, we have $\nu(B)=0$.
\end{remark}

\subsection{Examples}

A wealth of interesting examples are also presented in the recent paper of Janson \cite{Jans}. Here we include a very short list of examples to serve as illustration. Some of these examples are taken from Janson's paper. All these examples are easy to work out, we leave the simple calculations to the reader. For more examples of limits of \textbf{random finite trees} using metric sampling see Section~9-14 of \cite{Jans}.

\begin{example}[Example 7.1 \cite{Jans}]
Let $P_n$ be the \textbf{path} on $n$ vertices. These paths converge in metric sampling to the measured real tree $(I,d,\mu)$, where $I=[0,1]$ is the unit interval
with the Euclidean metric $d$ and Lebesgue measure $\mu$ on the Borel sets of $I$. This can be realized as a dendron $(I,d,\nu)$, where $\nu$ is the push forward of the measure $\mu$ along the inclusion $\iota:I\to I\times[0,\infty)$ defined by $f(x)=(x,0)$. 
\end{example}
\begin{example}[Example 7.2 \cite{Jans}]
For the simplest example where the limit is \emph{not} a measured real tree, consider the sequence of \textbf{stars}, consisting of the $n$-vertex star $K_{1,n-1}$ for all $n$. These stars converge in metric sampling to the dendron $(\Upsilon,d,\nu)$, where $(\Upsilon,d)$ is the one point metric space with point $p$ and $\nu$ is concentrated on the single point $(p,1/2)$.
\end{example}
\begin{example}[Example 7.3 \cite{Jans}]
Let $B_n$ be the \textbf{binary tree} of height $n-1$ with $2^n-1$ vertices. These binary trees also converge in metric sampling to the same dendron $(\Upsilon,d,\nu)$ we saw in the previous example.
\end{example}
\begin{example}
We obtain a more complicated limit dendron if we replace the edges of the binary tree with paths of variable lengths. Let us obtain $B'_n$ by replacing each edge $e$ of the binary tree $B_n$ from the preceding example with a path of length $\lfloor n^2/k^2\rfloor$, where $e$ connects a vertex of distance $k$ from the root with one at distance $k-1$. These finite trees form a convergent sequence in metric sampling. The limit dendron is $(B,d,\nu)$, where $(B,d)$ is a real tree obtained from the infinite binary tree by making each edge $e$ of the combinatorial tree in distance $k-1$ from the root into an interval of length $C/k^2$ where $C=\frac12\sum_{i=1}^\infty1/i^2$. The distribution $\nu$ is concentrated on the set $H\times\{0\}$ for the Cantor set $H$ formed by the points of $B$ at distance $1/2$ from the root (the ``limit points'' of the infinite branches) and it is uniform there.
\end{example}
\begin{example}
A limit dendron with a ``more uniform'' distribution can be obtained from \textbf{combs}. Let $C_n$ be the $n^2$-vertex tree obtained from the $n$-vertex path $P_n$ by attaching pairwise disjoint $n$-vertex paths at every vertex. These combs form a convergent sequence in metric sampling. The limit dendron is $(I,d,\nu)$, where $(I,d)$ is the interval of length $1/3$ and $\nu$ is concentrated on the rectangle $I\times[0,1/3]$ and is uniform there.
\end{example}
\begin{example}
The underlying real trees in the limit dendrons of the previous examples are all compact. Here is a simple way to get a non-compact underlying tree. Let $D_n$ be the depth two finite tree where the root has $n$ children: $P_1,\ldots P_n$ and $P_i$ has $2^i$ children for $i=1,\dots,n$. The trees $D_n$ also form a convergent sequence in metric sampling with the limit dendron being $(S,d,\nu)$, where the real tree $(S,d)$ is a star consisting of the infinitely many intervals $(r,p_i)$, each of length $1/4$. The distribution $\nu$ is concentrated on the points $(p_i,1/4)$ with $\nu((p_i,1/4))=2^{-i}$ for $i=1,2\dots$.
\end{example}
\medskip

In Section~\ref{semi} we introduce semi-measured real trees, a technical relaxation of measured real trees. In Sections~\ref{metult} and \ref{mesult} we recall the metric ultraproduct and the ultraproduct of measure spaces, respectively, especially as they apply to semi-measured real trees. We prove Theorems~\ref{tetel1}, \ref{tetel2} and \ref{unique} in Sections~\ref{pr1}, \ref{pr2} and \ref{un}, respectively.

\section{Semi-measured real trees}\label{semi}

\begin{definition}
We denote the open ball of radius $r$ around a point $x$ in a metric space by $B_x(r)$.

We call the triple $(T,d,\mu)$ a {\bf semi-measured metric space} if $(T,d)$ is a metric space, $(T,\mu)$ is a probability measure space and all the balls in $(T,d)$ are $\mu$-measurable. A semi-measured metric space $(T,d,\mu)$ is a {\bf quasi metric measure space} (or {\bf quasi mms} for short) if $d$ is $\mu^2$-measurable.
\end{definition}

\begin{remark} Note that both semi-measured metric spaces and a quasi metric measure spaces are relaxations of the well established notion of metric measure spaces. For metric measure spaces one first requires that the underlying metric space is separable and complete, second that the distribution is a probability Borel measure on this metric space. Here we dropped the first requirement and relaxed the second, allowing some Borel sets to be non-measurable. The requirement for the balls to be measurable in a semi-measured metric space is equivalent to requiring that the single variable distance function $d_x(y)=d(x,y)$ is $\mu$-measurable for all points $x\in T$ and as such, it is weaker even than requiring that the bivariate distance function $d$ is $\mu^2$-measurable as needed for a quasi mms.

Note also, that the sampling measures $\tau_r$ are defined for a quasi mms but not (in general) for semi-metric measure spaces.
\end{remark}

\begin{definition}
A {\bf semi-measured real tree} is a semi-measured metric space $(T,d,\mu)$, where $(T,d)$ is an $\R$-tree. Similarly, a {\bf quasi measured real tree} is a quasi mms $(T,d,\mu)$, where $(T,d)$ is a real tree.

If $Y$ is non-empty, closed, connected subset of $T$ in a real tree $(T,d)$ we define the retraction
$\pi_Y=\pi^{T,d}_Y:T\to Y$ by setting $\pi_Y(t)$ be the unique closest point to $t$ in $Y$ (the existence of which is stated in part~(2) of the following lemma).
\end{definition}

\begin{lemma}\label{subtree}
Let $(T,d)$ be a real tree and $Y$ be a non-empty, closed, connected subset of $T$.  
\begin{enumerate}
\item
$Y$ with the restriction of $d$ is a real tree.

\item
The retraction $\pi_Y$ is well defined and we have $\pi_Y(t)\in[t,y]$ for all  $t\in T$ and $y\in Y$.

\item
If $B$ is a connected component of $T\setminus Y$ and $x\in B$, then $B$ is a $\pi_Y(x)$-branch of $T$.
 
\item
We have $d(x,y)=d(x,\pi_Y(x))+d(\pi_Y(x),\pi_Y(y))+d(\pi_Y(y),y)$ for all $x,y\in T$ unless $x$ and $y$ are in the same connected component of $T\setminus Y$. In particular, we have $d(\pi_Y(x),\pi_Y(y))\le d(x,y)$ for all $x,y\in T$ making $\pi_Y$ continuous.

\item
Any branch $B$ of $T$ is in the $\sigma$-algebra generated by the balls of $T$, so it is $\mu$-measurable for any semi-measured real tree $(T,d,\mu)$.
\end{enumerate}
\end{lemma}

\proof
A closed subspace $Y$ of $T$ is a complete metric space. If $[x,y]\subseteq Y$ for all $x,y\in Y$ then $Y$ is a real tree, otherwise $Y$ is not connected. This proves part~(1).

Let us fix $t\in T$ and for a point $y\in Y$ define $p_y$ to be the unique closest point to $t$ in $Y\cap[y,t]$. As the segment $[y,t]$ is isometric to an interval and $Y$ is closed, this exists. If $p_{y_1}\ne p_{y_2}$ for $y_1,y_2\in Y$, then $p_{y_1}$ and $p_{y_2}$ can be connected inside $Y$ (as $Y$ is connected) and also outside (through $t$), a contradiction. Thus, all the points $p_y$ coincide defining $\pi_Y(t)$ and proving part~(2).

For part~(3) consider the $\pi_Y(x)$-branch containing $x$. By part~(2) it is disjoint from $Y$, hence it is $B$.

For $x,y\in T$, the union of the segments $[x,\pi_Y(x)]$, $[\pi_Y(x),\pi_Y(y)]$ and $[\pi_Y(y),y]$ connect $x$ to $y$. If some two of these three segments intersect in more than their end points, then it must be the first and last ones (as the middle segment is contained in $Y$), and then $x$ and $y$ are in the same connected component of $T\setminus Y$. If no such non-trivial intersection occurs, then the union of the three segments is homeomorphic to an interval, so the union must be $[x,y]$ itself, proving the formula for $d(x,y)$ in part~(4). If the formula applies, it implies the bound $d(\pi_Y(x),\pi_Y(y))\le d(x,y)$. Otherwise $\pi_Y(x)=\pi_Y(y)$ by part~(3) and the bound holds again.

Let $B$ be a $p$-branch of the real tree $(T,d)$ and $x\in B$. Consider a sequence of intermediate points $p_n$ of $[x,p]$ tending to $p$. It is easy to see that $B=\{y\in T\mid\exists n:d(y,p_n)<d(y,p)\}$. This makes $B=\bigcup_{n,r}(B_{p_n}(r)\setminus B_p(r))$, where the union is taken for all $n$ and all rational numbers $r>0$. The formula proves part~(5). \qed

The goal of this section is to study semi-measured real trees $(T,d,\mu)$ and to define the associated long dendrons for them. Note that the associated long dendron is, in fact, a dendron whenever the essential diameter of $(T,d,\mu)$ is at most one, that is, when $\mu(\{x\in T\mid\mu(\{y\in T\mid d(x,y)>1\})>0\})=0$.

\begin{definition}
Let $(T,d,\mu)$ be a semi-measured real tree. The point $p\in T$ is an {\bf inner point} of $T$ if there is no $p$-branch $B$ of
$T$ such that $\mu(B)=1$. The {\bf core} of $T$, $\core(T)$ is the closure
of the set of inner points of $T$. We write $\pi^T$ as a shorthand for the retraction $\pi_{\core(T)}$.

For a semi-measured real tree we define the {\bf associated long dendron} $D$ and {\bf associated projection} $\alpha:T\to A_D$ as follows. For $p\in T$ we set $\alpha(p)=(\pi^T(p),d(p,\pi^T(p)))$. We set $D=(\core(T),d,\nu)$. Here we slightly abuse notation by denoting the restriction of $d$ to the $\core(T)$ by $d$ again. We define the Borel probability measure $\nu$ as the push-forward of the measure $\mu$ along the map $\alpha$.
\end{definition}

The definition of inner point makes sense because the branches are measurable by Lemma~\ref{subtree}(5).

In order to show that the definition of the associated long
dendron makes sense we need to prove a series of lemmas.
In Lemmas~\ref{core1} and \ref{connected} we will show that $\core(T)$ is non-empty and connected, respectively. As the core is closed by definition it is a real tree by Lemma~\ref{subtree}(1) and so the retraction $\pi^T$ exists and the associated projection $\alpha$ is also defined. For the definition of the measure $\nu$ in the associated long dendron to make sense we further need that the associated projection $\alpha$ is measurable if considered as a map from $(T,\mu)$ to the Borel space on $A_D$. This is stated in Lemma~\ref{alphameasurable} below. Finally, we show that the associated long dendron is indeed a long dendron (that is, it satisfies the positivity condition) in Lemma~\ref{pos} below.

\begin{lemma} \label{core1}
The set $\core(T)$ is non-empty. Furthermore, a branch $B$ of $T$ intersects
$\core(T)$ if and only if $\mu(B)>0$.
\end{lemma}

\proof
It is enough to prove the second statement as if there are no positive measure branches, then all points of $T$ must be inner points. (This scenario is only possible if $T$ has a single point only and thus it has no branches at all.)

Let $p\in T$ and let $B$ be a $p$-branch. For any $q\in B$
the $q$-branch $B'$ containing $p$ satisfies $B'\supseteq T\setminus B$. Thus, if $\mu(B)=0$, then $\mu(B')=1$ and no point in $B$ is an inner point. As $B$ is open, this implies that $\core(T)$ is
disjoint from $B$ as claimed. Now, let us assume that $\mu(B)>0$. We will proceed by contradiction. Suppose that $B$
contains no inner points. For $t\in B$ let $B_t$ be the unique $t$-branch of measure $1$, and let $C_t$ be the $t$-branch containing $p$. Let
$S=\{t\in B\mid B_t\ne C_t\}$. Note that $\mu(C_t)=1$ if $t\in
B\setminus S$ and $\mu(C_t)=0$ if  $t\in S$.
\noindent
Let $z=\sup\{d(p,t)\mid t\in S\cup\{p\}\}$. We claim that $z$ is a maximum,
that is, there exists $c\in S\cup\{p\}$ with $d(p,c)=z$. 
Otherwise we would have $c_i\in S$ such
that $d(p,c_i)$ tend to $z>0$ as $i$ tends to infinity. We have
$\mu(\bigcap^\infty_{i=1}B_{c_i})=1$. Let $t\in\bigcap^\infty_{i=1}B_{c_i}$. All the points $c_i$
are in the segment $[p,t]$, so they must tend to a point $c\in[p,t]$ with
$d(p,c)=z$. We have $C_c=\bigcup^\infty_{i=1}C_{c_i}$ so we must have $\mu(C_c)=0$ and
$c\in S$ as claimed.
\noindent
Now, let $B'=B$ if $c=p$ (that is $S=\emptyset$) and $B'=B_c$ otherwise. Take an
arbitrary point $y\in B'$ and a sequence of intermediate points $y_i$ of the segment
$[c,y]$ tending to $c$. Since $d(p,y_i)> d(p,c)=z$, we have that $\mu(C_{y_i})=1.$ Notice that $\bigcap^\infty_{i=1}C_{y_i}$ is disjoint from
the positive measure set $B'$, therefore we cannot have $\mu(C_{y_i})=1$ for
all $i$, leading to a contradiction \qed

\begin{lemma}\label{connected}
The set $\core(T)$ is connected. Moreover, all intermediate points of the segment $[p,q]$ are inner points if $p,q\in\core(T)$.
\end{lemma}

\proof
Suppose now that $x$ is an intermediate point in $[p,q]$ for some $p,q\in\core(T)$. The points $p$
and $q$ lie in distinct $x$-branches and these $x$-branches have positive
measure by Lemma~\ref{core1}. Thus, $x$ must be an inner point as
stated. \qed

Let us fix a point $p_0\in\core(T)$ and let $Q$ be the set of inner points $q$
such that $d(p_0,q)$ is a rational number. 

\begin{lemma}\label{separable}
The real tree $\core(T)$ is separable. In particular, $Q$ is a countable dense
set in $\core(T)$ and furthermore $Q$ is dense in any proper segment in $\core(T)$.
\end{lemma}

\proof
For $p_0\ne q\in Q$, let $B_q$ be the $q$-branch containing $p_0$ and $C_q=T\setminus
B_q$. As $q$ is inner we have $\mu(C_q)>0$. Notice that for $p,q\in Q$ distinct with
$d(p_0,p)=d(p_0,q)$ the sets $C_p$ and $C_q$ are disjoint. Therefore, $Q$
contains a countable number of points in any fixed distance from $p_0$ and
must itself be a countable set. Consider a pair of distinct points $p,q\in\core(T)$ and let
$t=\pi_{[p,q]}(p_0)$. All intermediate points $s\in[p,q]$ are inner points by
Lemma~\ref{connected}. By Lemma~\ref{subtree} we have $d(p_0,s)=d(p_0,t)+d(t,s)$, and thus
$Q\cap[p,q]$ is dense in $[p,q]$ as claimed.

The fact that $Q$ is dense in any proper segment in $\core(T)$ implies that $Q$ is dense in $\core(T)$ itself unless $\core(T)$ is a singleton set. In this latter case $Q=\core(T)$. \qed

\begin{remark}
The same argument (with ``inner'' replaced with ``intermediate'') shows that every separable real tree has a countable subset that is dense in every proper segment.
\end{remark}

\begin{corollary}\label{corsep}
If $D=(T,d,\nu)$ is a long dendron, then the real tree $(T,d)$ is separable.
\end{corollary}

\proof
Let $\mu$ be push forward of the measure $\nu$ along the projection from $A_D$ to $T$. This makes $(T,d,\mu)$ a semi-measured real tree. By the definition of the long dendron, every branch of this real tree has positive measure. By Lemma~\ref{core1}, this implies that every branch intersects the core. But as the core is non-empty (by Lemma~\ref{core1} again) and closed (by definition), this implies that the core is $T$ itself. Now Lemma~\ref{separable} proves the claim of the corollary. \qed

\begin{lemma}\label{alphameasurable}
The associated projection $\alpha:T\to A_D$ is measurable.
\end{lemma}

\proof
Recall that for $p\in T$ we have $\alpha(p)=(\pi^T(p),d(\pi^T(p),p))$. We prove that both coordinate functions are measurable. Here $d(\pi^T(p),p)=\min_{q\in\core(T)}d(p,q)=\inf_{q\in Q}d(p,q)$ by Lemma~\ref{separable}. The function $f_q(p)=d(p,q)$ is measurable by definition, therefore so is $d(\pi^T(p),p)$.

We further need that the map $\pi^T$ considered from $(T,\mu)$ to the Borel space on $\core(T)$ is measurable. As $\core(T)$ is separable (Lemma~\ref{separable}) it is enough to check that the inverse image of the complement of a closed ball, namely the set $H=\{x\in T\mid d(\pi^T(x),p)>r\}$ is measurable for each $p\in\core(T)$ and $r>0$. For $x\in H$ there exists $q\in Q\cap[p,\pi^T(x)]$ with $d(q,p)>r$ by Lemma~\ref{separable}. This makes $d(x,p)>d(x,q)+r$. Clearly, for points $x$ outside $H$ and $q\in Q$ we have $d(x,p)\le d(x,\pi^T(x))+d(\pi^T(x),p)\le d(x,q)+r$. This proves $H=\bigcup_{q\in Q}\{x\in T\mid d(x,p)>d(x,q)+r\}=\bigcup_{q\in Q,s}(B_q(s)\setminus B_p(r+s))$, where $s$ ranges over the positive rational numbers in the last expression. As the balls are $\mu$-measurable by definition, so is $H$. \qed

\begin{lemma}\label{pos}
The associated long dendron $D=(\core(T),d,\nu)$ is indeed a long dendron.
\end{lemma}

\proof
We have already seen that $(\core(T),d)$ is a real tree and $\nu$ is a Borel probability measure on $A_D$. It remains to show that for any branch $B$ of $\core(T)$ we have $\nu(B\times[0,\infty))>0$. By the definitions of $\nu$ and $\alpha$ we have $\nu(B\times[0,\infty))=\mu(\{x\in T\mid\pi^T(x)\in B\})$.  Let $B$ be a $p$-branch of $\core(T)$, then $\{x\in T\mid\pi^T(x)\in B\})$ is the $p$-branch $B'$ of $T$ containing $B$. $B'$ intersects $\core(T)$, so by Lemma~\ref{core1} we have $\mu(B')>0$. \qed

We call the connected components of $T\setminus\core(T)$ the {\bf feathers} of $T$.

\begin{lemma}\label{feather}
Let $(T,d,\mu)$ be a semi-measured real tree with associated long dendron $D=(\core(T),d,\nu)$ and associated projection $\alpha$. 
The feathers of $T$ are exactly the measure zero $p$-branches for points $p\in\core(T)$. We have $d(x,y)=d_D(\alpha(x),\alpha(y))$ for $x,y\in T$ unless $x$ and $y$ are in the same feather of $T$.
\end{lemma}

\proof
Both statements follow directly from Lemmas~\ref{subtree} and \ref{core1}. \qed

A {\bf measured real tree} is a metric measure space $(T,d,\mu)$ where $(T,d)$ is a real tree. As we saw, a measured real tree is also a semi-measured real tree, so the associated long dendron and the associated projection are defined.

For a function $\alpha:A\to B$ and $r\in\N$ we write $\alpha^r$ for the function that acts coordinate-wise on $A^r$, that is $\alpha^r:A^r\to B^r$, $\alpha^r(x_1,\dots,x_r)=(\alpha(x_1),\dots,\alpha(x_r))$. 

\begin{lemma}\label{measureequal}
For a quasi mms $(T,d,\mu)$ and associated long dendron $D$ we have $\tau_r(T,d,\mu)=\tau_r(D)$ for all $r\ge1$.
\end{lemma}

\proof
Let $D=(\core(T),d,\nu)$ and let $\alpha:T\to A_D$ be the associated projection. By Lemma~\ref{feather} we have $d(x,y)=d_D(\alpha(x),\alpha(y))$ for $x,y\in T$ unless $x$ and $y$ are in the same feather of $T$. Also by the same lemma the feathers have zero measure, so we have $\mu(\{y\in T\mid d(x,y)\ne d_D(\alpha(x),\alpha(y))\})=0$ for all $x\in T$ and hence $\mu^2(\{(x,y)\in T^2\mid d(x,y)\ne d_D(\alpha(x),\alpha(y))\})=0$. As a consequence, we have $\mu^r(X)=0$ for the set $X=\{x\in T^r\mid\rho_r^{T,d}(x)\ne\rho_r^{A_D,d_D}(\alpha^r(x))\}$.

For a Borel set $H\subseteq M_r$ we have $\tau_r(T,d,\mu)(H)=\mu^r(\{x\in T^r\mid\rho_r^{T,d}(x)\in H\})$ and $\tau_r(D)(H)=\nu^r(\{x\in A_D^r\mid\rho_r^{A_D,d_D}(x)\in H\})=\mu^r(\{x\in T^r\mid\rho_r^{A_D,d_D}(\alpha^r(x))\in H\})$. As the symmetric difference of the sets $\{x\in T^r\mid\rho_r^{T,d}(x)\in H\}$ and $\{x\in T^r\mid\rho_r^{A_D,d_D}(\alpha^r(x))\in H\}$ is contained in the zero measure set $X$, their $\mu^r$ measures agree, so $\tau_r(T,d,\mu)(H)=\tau_r(D)(H)$ as needed. \qed

\begin{remark}\label{meareal}
Given a (long) dendron $D$ it is natural to look for a quasi measured real tree $T_D$ with $D$ as its associated long dendron. By Lemma~\ref{measureequal} we could use $T_D$ in place of $D$ for the sampling limit in Theorems~\ref{tetel1} and \ref{tetel2}. If the construction of $T_D$ is canonical we can also ensure the limit is unique as in Theorem~\ref{unique}.

Let $D=(T,d,\nu)$ be a long dendron. Recall that $d_D:A_D^2\to\R$ is not a metric as $d_D(x,x)$ is often positive. As a first attempt to fix this one could consider the distance function $d'$ as defined in Equation~(\ref{d'}) in Section~\ref{un}. 
This is indeed a distance and $(A_D,d')$ is a real tree. It is not separable (unless $T$ is trivial), so $(A_D,d',\nu)$ is typically not a metric measure space, but it is a quasi measured real tree so both $\tau_r(D)$ and $\tau_r(A_D,d',\nu)$ are well defined. If $\nu(\{u\}\times(0,\infty))=0$ for all $u\in T$, then $d'$ and $d_D$ differ in a measure zero subset of $A_D^2$. In this case we have $\tau_r(D)=\tau_r(A_D,d',\nu)$ and further the associated long dendron for the measured real tree $(A_D,d',\nu)$ is $D$ and the associated projection is the identity on $A_D$. Note however, that $\nu$ is a Borel measure on $A_D$ with its product structure, but unless $|T|=1$ it is \emph{not} a Borel measure on the real tree $(A_D,d')$ as the latter space is non-separable with many not $\nu$-measurable open sets.

The simple approach above does not work if $\nu(\{u\}\times(0,\infty))>0$ for some points $u\in T$. The following, slightly more involved method always works. We define the quasi measured real tree $T_D=(T^*,d^*,\mu)$ as follows.
$$T^*=\{(x,y,z)\mid x\in T,\,(y,z)\in\{(0,0)\}\cup((0,1)\times(0,\infty))\}$$
$$d^*((x_1,y_1,z_1),(x_2,y_2,z_2))=\left\{\begin{array}{ll}|z_1-z_2|&\hbox{ if }(x_1,y_1)=(x_2,y_2)\\d(x_1,x_2)+z_1+z_2&\hbox{ otherwise.}\end{array}\right.$$
We define $\mu$ as a Borel probability measure on $T^*$ considered as subset of the product $T\times[0,1)\times[0,\infty)$. As above, this is \emph{not} a Borel measure on the real tree $(T^*,d^*)$ because many open sets in the real tree are not measurable. For $H\subseteq T^*$ which is Borel in the former sense we define $\mu(H)$ through the function $f_H:A_D\to[0,1]$ as follows. We denote the Lebesgue measure on the reals by $\lambda$.
$$f_H(x,z)=\left\{\begin{array}{ll}0&\hbox{ if }z=0\hbox{ and }(x,0,0)\notin H\\1&\hbox{ if }z=0\hbox{ and }(x,0,0)\in H\\\lambda(\{y\in(0,1)\mid(x,y,z)\in H\})&\hbox{ if }z>0\end{array}\right.$$
$$\mu(H)=\int_{A_D}f_H\,d\nu$$
Now $T_D$ is a quasi measured real tree which has $D$ as its associated long dendron and $\alpha:T^*\to A_D$ given by $\alpha(x,y,z)=(x,z)$ as its associated projection. In particular, we have $\tau_r(T_D)=\tau_r(D)$ for all $r$.
\end{remark}

\section{The ultraproduct and the metric ultraproduct}\label{metult}

Let us recall the definitions of the ultralimit and the ultraproduct. We will use both the set theoretical and metric ultraproducts as well as the ultraproduct of probability measures (see next section). When referring to ultralimits or ultraproducts we will always use the same fixed nontrivial ultrafilter $\omega$ on the set $\N$ of positive integers. We use {\bf$\omega$-few} to refer to any subset of $\N$ that is not in $\omega$. We allow sequences $(x_n)_{n\in\N}$ to be undefined for $\omega$-few indices $n$.

For any Hausdorff space $H$ and points $x_n\in H$ for $n\in\N$ we define the ultralimit $\limo x_n$ to be $x\in H$ such that $x_n$ is outside any fixed neighborhood of $x$ for only $\omega$-few indices $n$. If the ultralimit exists it is clearly unique, and in case $H$ is compact, it does exist for every sequence $x_n$. Also, if $\lim_{n\to\infty}x_n$ exists, then so does $\limo x_n$ and they agree. Further, if $\limo x_n$ exists and $g:H\to H'$ is a continuous function to another Hausdorff space $H'$, then $\limo g(x_n)=g(\limo x_n)$.

The {\bf set-theoretical ultraproduct} $\uT=\prodo T_n$ of the sequence of sets $(T_n)_{n\in\N}$
consists of the equivalence classes of sequences
$(x_n)_{n\in\N}$, $x_n\in T_n$, where $(x_n)_{n\in\N}$ and $(y_n)_{n\in\N}$ are equivalent if $x_n\ne y_n$ for $\omega$-few indices $n$. We denote the class of the sequence $(x_n)_{n\in\N}$ by $[(x_n)_{n\in\N}]$.

Let $A_n\subseteq T_n$ for each $n$. Clearly, each element $\ux$ of $\uA=\prod_\omega A_n$ (as an equivalence class) is contained in a distinct element of $\uT=\prodo T_n$. We identify $\ux$ with the element of $\uT$ containing it making $\uA\subseteq\uT$.

For sets $S_n$, $T_n$ we identify the ultraproduct of their direct products $\prodo(S_n\times T_n)$ with the direct product $\prodo S_n\times\prodo T_n$ by identifying $[(x_n,y_n)_{n\in\N}]$ in the former space with $([(x_n)_{n\in\N}],[(y_n)_{n\in\N}])$ in the latter. We make the same identification in ultraproducts of products of finitely many factors. In particular, we identify $\prodo(S_n^k)$ with $(\prodo S_n)^k$ for $k\ge2$.

Let $H$ be a Hausdorff space and $f_n:T_n\to H$ be arbitrary functions. We write $\limo f_n$ for the function $\uf$ defined by $\uf([(x_n)_{n\in\N}])=\limo f_n(x_n)$. The ultralimit does not depend on the choice of the sequence $(x_n)_{n\in\N}$ representing the class $[(x_n)_{n\in\N}]$, but it may be undefined for some classes if $H$ is not compact. In this case $\uf$ is only defined on a subset of $\prodo T_n$.

The metric ultraproduct of metric spaces $(T_n,d_n)$ was introduced by van den Dries and Wilkie in \cite{DW}. First we define $\ud=\limo d_n$. For this we consider the distance functions $d_n$ with values in the compact space $\RR=\R\cup\{-\infty,\infty\}$. This makes $\ud:\prodo(T_n^2)\to\RR$. We identified $\prodo(T_n^2)$ with $\uT^2$, where $\uT=\prod_\omega T_n$, so $\ud:\uT^2\to\RR$. Clearly, $\ud$ is symmetric, non-negative and satisfies the triangle inequality, that is, $\ud$ is a {\bf pseudometric} possibly containing zero and infinite distances. The {\bf metric ultraproduct} $\bT=\hprodo T_n$ is obtained from $\uT$ by factoring out the equivalence relation of having zero distance. Let $\hpi:\uT\to\bT$ be the {\bf natural projection} mapping a point in $\uT$ to its equivalence class in $\bT$. We write $\widehat{[(p_n)_{n\in\N}]}$ to denote $\hpi([(p_n)_{n\in\N}])$. Clearly, the pseudodistance $\ud$ defines a pseudodistance $\bd$ on $\bT$ by setting $\bd(\hpi(\ux),\hpi(\uy))=\ud(\ux,\uy)$. Here $\bd$ is a pseudometric in which distinct points have a positive distance but infinite distances may still appear. Clearly,  $\diam(\bT,\bd)=\limo\diam(T_n,d_n)$. If this value is finite, we call the sequence $(T_n,d_n)_{n\in\N}$ {\bf uniformly bounded}. In this case $(\bT,\bd)$ is a metric space (no infinite distances). In general, ``having finite distance'' is an equivalence relation on $\bT$ and each equivalence class $\bX$ is made into a metric space by the restriction of $\bd$. We call these metric spaces the {\bf clusters} of the metric ultraproduct $\bT$.

Let $A_n\subseteq T_n$ for each $n$. As we identified $\prodo A_n$ with a subset of $\uT$ we similarly identify $\hprodo A_n$ with the corresponding subset of $\bT$.

Note that the ultraproduct objects will always  be typeset in bold and we will put the metric
ultraproduct objects under the $\hat{}$ sign. Recall that we use normalized distance on finite trees and all those distances are bounded by $1$. Therefore, sequences of finite trees are uniformly bounded. We could have restricted attention to uniformly bounded sequences of metric spaces in this paper to avoid dealing with several clusters. We allow unbounded diameter for more generality and because having several clusters does not significantly increase complexity.

\begin{lemma}\label{treeprod}
Any cluster $\bX$ of the metric ultraproduct $(\bT,\bd)$ of real trees $(T_n,d_n)$ is a real tree. For two points $\bx=\widehat{[(x_n)_{n\in\N}]}$ and $\by=\widehat{[(y_n)_{n\in\N}]}$ in $\bX$ we have $[\bx,\by]=\hprodo[x_n,y_n]$.
\end{lemma}

\proof As $\bX$ is obtained as a metric ultraproduct, it is complete, see \cite{DW}. For $\bx,\by\in\bX$ the segments $[x_n,y_n]$ in $T_n$ are isometric to a real interval of length $d_n(x_n,y_n)$ and we have $\limo d_n(x_n,y_n)=\bd(\bx,\by)<\infty$. Therefore, the metric ultraproduct $\bl=\hprodo[x_n,y_n]$ of these segments is isometric to a real interval of length $\bd(\bx,\by)$. It is contained in $\bX$ and contains $\bx$ and $\by$, so all we need to establish to prove the lemma is that all the intermediate points of $\bl$ separate $\bx$ from $\by$.

So let $\bp=\widehat{[(p_n)_{n\in\N}]}\in\bl\setminus\{\bx,\by\}$ with $p_n\in[x_n,y_n]$. Let $A_n$ be the $p_n$-branch in $T_n$ containing $x_n$ and $B_n=T_n\setminus A_n$. (For the $\omega$-few indices $n$ where $p_n=x_n$ the sets $A_n$ and $B_n$ are not defined.) Let $\bA=\hprodo A_n$ and $\bB=\hprodo B_n$. We have $\bA\cup\bB=\bT$ as for any $\bz=\widehat{[(z_n)_{n\in\N}]}\in\bT$ we have either $z_n\notin A_n$ for $\omega$-few indices making $\bz\in\bA$ or $z_n\notin B_n$ for $\omega$-few indices making $\bz\in\bB$. Take any $\ba\in\bA$ and $\bb\in\bB$. From the similar equations in $T_n$ we have $\bd(\ba,\bb)=\bd(\ba,\bp)+\bd(\bp,\bb)$. This means that any point $\ba\in\bA\setminus\{\bp\}$ is separated from $\bB$ by the positive distance of $\bd(\ba,\bp)$ and similarly any point $\bb\in\bB\setminus\{\bp\}$ is separated from $\bA$ by the positive distance of $\bd(\bp,\bb)$. Therefore the sets $(\bA\cap\bX)\setminus\{\bp\}$ and $(\bB\cap\bX)\setminus\{\bp\}$ are disjoint and open and form a partition of $\bX\setminus\{\bp\}$. The point $\bx$ is in the former set, while $\by$ is in the latter, so $\bp$ separates them in $\bX$ as claimed. \qed

\section{The ultraproduct of probability spaces}\label{mesult}

Let us recall the ultraproduct of probability measures. Let $\mu_n$ be a probability measure on a $\sigma$-algebra $\cA_n$ over a set $T_n$ for each $n\in\N$. The ultraproduct sets $\prodo A_n$ with $A_n\in\cA_n$ form a Boolean algebra $\cP$ on $\uT=\prodo T_n$.
Additionally, we have a finitely additive measure $\mu_{\cP}$
on $\cP$ given by $\mu_{\cP}(\prodo A_n)=\limo\mu_n(A_n)$. This finitely additive measure can be extended to 
a $\sigma$-algebra containing $\cP$ (see \cite{ESZ}) as follows.
We call $\uN\subset\uT$ a {\bf nullset} if for any $\e>0$ there exists
an element $\uA\in\cP$ such that $\uN\subseteq\uA$ and $\mu_{\cP}(\uA)<\e$.
A set $\uM\subset\uT$ is called {\bf measurable} if there exists $\uP\in\cP$ such
that the symmetric difference $\uP\triangle\uM$ is a nullset. The family of measurable sets $\cM_{\uT}$
form a $\sigma$-algebra with a probability measure $\mu=\prodo\mu_n$ (the ultraproduct of the measures $\mu_n$),  where we define $\mu(\uM)=\mu_{\cP}(\uP)$. Hence, we made the ultraproduct space $\uT$ into a probability measure space $(\uT,\mu)$. Note that the paper \cite{ESZ} only dealt with ultraproducts of finite spaces with uniform probability measures, so the following lemma is only stated there for this special case. However, the proof of the lemma does not use this restriction on the factors of the ultraproduct and the same proof applies verbatim for the more general case stated here.

\begin{lemma} [Lemma 5.1 \cite{ESZ}]\label{functionlimit}
Let $(T_n,\mu_n)$ be probability measure space and $f_n:T_n\to\RR$ be a $\mu_n$-measurable function for $n\in\N$. Let $\uT=\prodo T_n$, $\mu=\prodo\mu_n$. In this case the function $\uf=\limo f_n:\uT\to\RR$ is $\mu$-measurable. If $(f_n)_{n\in\N}$ is uniformly bounded (that is $\sup_{n,x}|f_n(x)|<\infty$), then we also have
$$\limo\int_{T_n}f_n\,d\mu_n=\int_{\uT}\uf\,d\mu\,.$$
\end{lemma}

Let $(S_n,\mu_n)$ and $(T_n,\nu_n)$ be probability measure spaces for each $n\in\N$. Let $\uS=\prodo S_n$, $\uT=\prodo T_n$. Recall that we have identified $\uS\times\uT$ with $\prodo(S_n\times T_n)$. We have two probability measures on this set. First we have $\mu_{\uS}\times\mu_{\uT}$, where $\mu_{\uS}=\prodo\mu_n$ and $\mu_{\uT}=\prodo\nu_n$. But we also have $\mu_{\uS\times\uT}=\prodo(\mu_n\times\nu_n)$. It is easy to see that the letter measure extends the former, that is, any $\mu_{\uS}\times\mu_{\uT}$-measurable set $H$ is also $\mu_{\uS\times\uT}$-measurable and we have $\mu_{\uS\times\uT}(H)=(\mu_{\uS}\times\mu_{\uT})(H)$.

In general, a $\mu_{\uS\times\uT}$-measurable set is not necessarily
$\mu_{\uS}\times\mu_{\uT}$-measurable, but we still have the following form of
Fubini's theorem (see \cite{Loeb} and also \cite{Evan})

\begin{lemma}\label{fubini}
In the setting above for any $\mu_{\uS\times\uT}$-measurable set $H$ the sections $H_{\ux}=\{\uy\mid(\ux,\uy)\in H\}$ are $\mu_{\uT}$-measurable unless $\ux\in\uZ$ for some $\uZ\subseteq\uS$ with $\mu_{\uS}(\uZ)=0$ and we have
$$\mu_{\uS\times\uT}(H)=\int_{\uS\setminus\uZ}\mu_{\uT}(H_{\ux})\,d\mu_{\uS}(x).$$
\end{lemma}

Let $(T_n,d_n,\mu_n)$ be semi-measured metric spaces. Let $\uT=\prodo T_n$ be the set theoretic ultraproduct, let $(\bT,\bd)$ be the metric ultraproduct of the metric spaces $(T_n,d_n)$ with $\hpi:\uT\to\bT$ being the natural projection.

Let $\mu$ be the ultraproduct of the measures $\mu_n$. We call the sets $H\subseteq\bT$ {\bf proper} if $\hpi^{-1}(H)$ is $\mu$-measurable and define $\hmu(H)=\mu(\hpi^{-1}(H))$ in this case. Clearly, the proper sets form a $\sigma$-algebra on $\bT$ and $\hmu$ is a probability measure on it, the push-forward of $\mu$ along $\hpi$.

Note that if $S_n\subseteq T_n$ is $\mu_n$-measurable for $n\in\N$ and $\bS=\hprodo S_n\subseteq \bT$ is proper, then we have $\hmu(\bS)\ge\limo\mu_n(S_n)$. Indeed, $\limo\mu_n(S_n)=\mu(\uS)$, where $\uS=\prodo S_n$ while $\hmu(\bS)=\mu(\{\up\in\uT\mid\hpi(\up)\in\hpi(\uS)\})$ and the set $\{\up\in\uT\mid\hpi(\up)\in\hpi(\uS)\}$ clearly contains $\uS$.

Let $\ux=[(x_n)_{n\in\N}]\in\uT$. By definition, all the single variable distance functions $f_n:T_n\to\R$ given by $f_n(y)=d_n(x_n,y)$ are $\mu_n$-measurable. By Lemma~\ref{functionlimit} the function $\uf=\limo f_n$ is $\mu$-measurable. Here $\uf:\uT\to\RR$ is the function $\uf(\uy)=\ud(\ux,\uy)$. Hence its factor $\bbf:\bT\to\RR$ given by $\bbf(\by)=\bd(\hpi(\ux),\by)$ is $\hmu$-measurable.

This, in particular, means that each cluster of $\bT$ is proper. We call a cluster {\bf dominant} if its measure is $1$. If a dominant cluster exists, then we say that the sequence $(T_n,d_n,\mu_n)$ is {\bf essentially bounded}. In this case we simply disregard the rest of the metric ultraproduct and consider $(\bX,\bd_{\bX},\hmu_{\bX})$ as the ultraproduct of the semi-measurable metric spaces $(T_n,d_n,\mu_n)$, where $\bX$ is the dominant cluster and $\bd_{\bX}$ and $\hmu_{\bX}$ are the restrictions of $\bd$ and $\hmu$, respectively. By the above observation the ultraproduct of essentially bounded semi-measured metric spaces is a semi-measured metric space. Note that if $(T_n,d_n)$ is uniformly bounded, then $(T_n,d_n,\mu_n)$ is essentially bounded, but the converse often fails. 

\begin{remark}\label{strong}
The ultraproduct of an essentially bounded sequence of metric measure spaces (or even of a uniformly bounded sequence of metric measure spaces) is typically not a metric measure space, not even a quasi mms. A simple example of this phenomenon is obtained from random graphs. Let us associate with the $n$ vertex graph $G_n$ the metric measure space $X_n$ on the $n$ vertices of $G_n$ with the uniform distribution, where the distance between adjacent vertices is $1$ and between non-adjacent vertices it is $2$. The diameter of these spaces is bounded by $2$, so the metric measure spaces $X_n$ form a uniformly bounded sequence. Their ultraproduct is clearly non-separable. Furthermore, with probability $1$ for uniform random graphs $G_n$, the ultraproduct is not even a quasi mms.

We introduced semi-measured metric spaces because the ultraproduct of an essentially bounded sequence of semi-measured metric spaces is always a semi-measured metric space. However, this detour is not necessary. It is not hard to see that the ultraproduct of an essentially bounded sequence of measured real trees (or even quasi measured real trees) is always a quasi measured real tree. We could simply use this ultraproduct as the sampling limit (as opposed to the associated long dendron) by Lemma~\ref{measureequal}. But this approach does not give a unique sampling limit and the non-separability also complicates things. If one still finds a quasi measured real tree as a more esthetically pleasing sampling limit, this can be achieved without losing uniqueness through the process explained in Remark~\ref{meareal}.
\end{remark}

\section{The proof of Theorem \ref{tetel1}}\label{pr1}

Let $T$ be a nontrivial finite tree. Consider the real tree $(T',d)$ obtained from $T$ by turning its edges into equal length real segments and letting $d$ be the normalized distance such that $\diam(T',d)=1$. Let $\mu$ be the Borel probability measure on $T'$ that is concentrated on the vertices of $T$ and is uniform there. Now $(T',d,\mu)$ is a measured real tree and we clearly have $\tau_r(T',d,\mu)=\tau_r(T)$ for each $r\ge1$.

Theorem~\ref{tetel1} is a consequence of Theorem~\ref{tetel1+} below. Indeed, if $(T_n)_{n\in\N}$ is a convergent sequence of finite trees, then the corresponding sequence $(T'_n,d_n,\mu_n)$ of measured real trees is uniformly bounded, so Theorem~\ref{tetel1+} finds their limit long dendron. This long dendron $D$ is, in fact, a dendron, since $\diam(T'_n,d_n,\mu_n)=1$ for all $n\in\N$, so $\tau_2(D)$ is concentrated on matrices with no entry exceeding $1$.

\begin{theorem}\label{tetel1+}
Let $(T_n,d_n,\mu_n)_{n\in\N}$ be an essentially bounded sequence of measured real trees. For any $r\in\N$ we have
$$\limo\tau_r(T_n,d_n,\mu_n)=\tau_r(D)\,,$$
where $D$ is the long dendron associated with the semi-measured real tree $(\bX,\bd,\hmu)$ obtained as the ultraproduct of the measured real trees $(T_n,d_n,\mu_n)$ and the ultralimit is understood in the space $\prob(M_r)$ with the weak topology.
\end{theorem}

\proof Let $\hpi:\uT\to\bT$ be the natural projection from the set theoretic ultraproduct $\uT=\prodo T_n$ to the metric ultraproduct $\bT=\hprodo T_n$. Note that $\bX$ is the dominant cluster of $\bT$. Let $D=(X^*,\bd,\nu)$, and let the associated projection be $\alpha:\bX\to A_D$. For the proof we need to show that for any bounded continuous function $g:M_r\to\R$ we have
\begin{equation}\label{cel}
\int_{M_r}g\,d\tau_r(D)=\limo\int_{M_r}g\,d\tau_r(T_n,d_n,\mu_n)\,.
\end{equation}

Let us set $\mu=\prodo\mu_n$ and $\uX=\hpi^{-1}(\bX)$. For the left hand side of Equation~(\ref{cel}) we use the fact that we defined several of our measures as push-forwards of other measures.
\begin{eqnarray}\label1
\int_{M_r}g(z)\,d\tau_r(D)(z)&=&\int_{A_D^r}g(\rho_r^{A_D,d_D}(x))\,d\nu^r(x)\\\label2
&=&\int_{\bX^r}g(\rho_r^{A_D,d_D}(\alpha(\bx)))\,d\hmu^r(\bx)\\\label3
&=&\int_{\uX^r}g(\rho_r^{A_D,d_D}(\alpha(\hpi(\ux))))\,d\mu^r(\ux)\,.
\end{eqnarray}
Namely, the equation in line~(\ref1) holds as $\tau_r(D)$ is defined as the push-forward of $\nu^r$ along $\rho_r^{A_D,d_D}$. Line~(\ref2) follows as $\nu$ is defined as a push-forward of $\hmu$ along $\alpha$. Finally, line~(\ref3) follows as $\hmu$ was itself defined as the push-forward of $\mu$ along $\hpi$. Here $x$ is an $r$-tuple from $A_D$, $\bx$ is an $r$-tuple from $\bX$ and $\ux$ is an $r$-tuple from $\uX$. The functions $\alpha$ and $\hpi$ act on $r$-tuples coordinate-wise.

For the right hand side of Equation~(\ref{cel}) we have a longer sequence of equations. We use some of the same notations as above and will explain each line separately together with any additional notation introduced there.
\begin{eqnarray}\label4
\limo\int_{M_r}g(z)d\tau_r(T_n,d_n,\mu_n)(z)&=&\limo\int_{T_n^r}g(\rho_r^{T_n,d_n}(x_n))\,d\mu_n^r(x_n)\\\label5
&=&\int_{\uT^r}\limo\,g(\rho_r^{T_n,d_n}(x_n))\,d\mu^{(r)}([(x_n)_{n\in\N}])\\\label6
&=&\int_{\uX^r}\limo\,g(\rho_r^{T_n,d_n}(x_n))\,d\mu^{(r)}([(x_n)_{n\in\N}])\\\label7
&=&\int_{\uX^r}g(\rho_r^{\uX,\ud}(\ux))\,d\mu^{(r)}(\ux)\\\label8
&=&\int_{\uX^r}g(\rho_r^{A_D,d_D}(\alpha(\hpi(\ux))))\,d\mu^{(r)}(\ux)
\end{eqnarray}

Line~(\ref4) follows again from the definition of $\tau_r(T_n,d_n,\mu_n)$ as push-forward measure.

To obtain line~(\ref5) we apply Lemma~\ref{functionlimit} to the measure spaces $(T_n^r,\mu_n^r)$. We have $\prodo T_n^r=\uT^r$. We wrote $\mu^{(r)}$ to denote $\prodo\mu_n^r$. Recall that this measure is an extension of $\mu^r$.

We have $\mu(\uX)=\hmu(\bX)=1$ and therefore $\uX^r$ is a full measure subset of $\uT^r$, so restricting the domain to it, as done in line~(\ref6), does not affect the integral.

Let $\ud$ be the restriction of $\limo d_n$ to $\uX^2$. Note that $\ud$ does not take infinite values. Let $\ux=(\ux^1,\dots,\ux^r)\in\uX^r$ with $\ux^i=[(x^i_n)_{n\in\N}]$ and $x^i_n\in T_n$. Let us write $x_n=(x^1_n,\ldots,x^r_n)$. We have
\begin{eqnarray*}
\limo\,\rho_r^{T_n,d_n}(x_n)&=&\limo(d_n((x_n)_i,(x_n)_j))_{1\le i,j\le r}\\
&=&(\limo\,d_n((x_n)_i,(x_n)_j))_{1\le i,j,\le r}\\
&=&(\ud(\ux_i,\ux_j))_{1\le i,j\le r}\\
&=&\rho_r^{\uX,\ud}(\ux)\,.
\end{eqnarray*}
As $g$ is continuous, this implies
$$g(\rho_r^{\uX,\ud}(\ux))=g(\limo\rho_r^{T_n,d_n}(x_n))=\limo\,g(\rho_r^{T_n,d_n}(x_n))\,,$$
so we have the exact same integral in line~(\ref7) as we had in line~(\ref6).

Consider a point $\ux\in\uX^r$ where the integrands of lines~(\ref7) and (\ref8) differ. Let $\bx=\hpi(\ux)\in\bX^r$. We must have $\rho_r^{\uX,\ud}(\ux)\ne\rho_r^{A_D,d_D}(\alpha^r(\bx))$. Both sides are matrices in $M_r$ with zero entries in the diagonal. The off-diagonal entry in position $(i,j)$ is $\ud(\ux_i,\ux_j)$ and $d_D(\alpha(\bx_i),\alpha(\bx_j))$, respectively. Consider the set $\uH_{i,j}\subseteq\uX^r$ consisting of points where these values differ and project it to coordinates $i$ and $j$. This projection from $(\uT^r,\mu^{(r)})$ to $(\uT^2,\mu^{(2)})$ is clearly measure preserving and the image is the zero measure set $\uH$ described in Lemma~\ref{nullmertek} below. Thus, $\mu^{(r)}(\uH_{i,j})=0$.

We saw that the integrands in lines~(\ref7) and (\ref8) agree outside the zero measure set $\bigcup_{i\ne j}\uH_{i,j}$, so the two integrals agree.

Finally, notice that we integrate the exact same function on the same domain in lines~(\ref3) and (\ref8). We use different measures but one of them is an extension of the other. As both integrals are defined, they must agree. This shows that Equation~(\ref{cel}) holds and finishes the proof of the theorem assuming Lemma~\ref{nullmertek} below. \qed

\begin{lemma}\label{nullmertek}
Using the notation of Theorem~\ref{tetel1+} and its proof we have $\mu^{(2)}(\uH)=0$ for the set $\uH=\{(\ux,\uy)\in\uX^2\mid\ud(\ux,\uy)\ne d_D(\alpha(\hpi(\ux)),\alpha(\hpi(\uy)))\}$.
\end{lemma}

\proof The right hand side of the inequality defining $\uH$ is $\mu^2$-measurable. The left hand side is $\ud=\limo d_n$, so it is $\mu^{(2)}$-measurable by Lemma~\ref{functionlimit}. This makes $\uH$ also $\mu^{(2)}$-measurable. Its slices are $\uH_\ux=\{\uy\in\uT\mid(\ux,\uy)\in\uH\}$. By Lemma~\ref{feather} and since $\ud(\ux,\uy)=\bd(\hpi(\ux),\hpi(\uy))$, we have $(\ux,\uy)\in\uH$ if and only if $\hpi(\ux)$ and $\hpi(\uy)$ are in the same feather of $(\bX,\bd,\hmu)$. This makes the slice $\uH_\ux$ either empty (if $\ux\in\uT\setminus\uX$ or $\hpi(\ux)\in\core(\bX)$) or $\hpi^{-1}(\bF)$ for a feather $\bF$. We have $\hmu(\bF)=0$ by Lemma~\ref{feather}, so $\mu(\uH_\ux)=0$ in both cases. We have $\mu^{(2)}(\uH)=\int_{\uT}\mu(\uH_\ux)\,d\mu(\ux)=\int_\uT0\,d\mu=0$ by Lemma~\ref{fubini}. \qed

\section{The proof of Theorem~\ref{unique}}\label{un}

We will use the following simple observations of compact subtrees exhausting some real trees.

\begin{lemma}\label{compact}
Let $(T,d,\mu)$ be a semi-measured real tree with $\mu(B)>0$ for all branches $B$ of $T$. Then there exists a countable set $Q\subseteq T$ that is dense in any proper segment of $T$. For all $\e>0$ there exists a compact subtree $Y$ of $T$ with $\mu(Y)>1-\e$.
\end{lemma}

\proof
Any intermediate point $p$ has more than one $p$-branches, each of positive measure, making $p$ an inner point. The intermediate points are dense in $T$ unless $|T|=1$, so we have $\core(T)=T$. This trivially holds even if $|T|=1$. The separability comes from Lemma~\ref{separable}.

Let us fix $p_0\in T$ and for a rational number $r>0$ let $Q_r$ be the set of inner points in $T$ of distance $r$ from $p_0$. Recall that the set $Q=\bigcup_rQ_r$ is a countable set that is dense in every proper segment of $T$ by Lemma~\ref{separable}.

To find the compact subtree with large measure let us fix constants $\e_r>0$ for the positive rational numbers $r$ and $\e'>0$ with $\e'+\sum_r\e_r<\e$. Let us find a closed ball $A$ around $p_0$ with $\mu(A)>1-\e'$.
Let us define $B_p$ for $p\in Q$ to be the $p$-branch of $T$ containing $p_0$ and $C_p=T\setminus(B_p\cup\{p\})$. Clearly, the open sets $C_p$ for $p\in Q_r$ are pairwise disjoint, so we have $\sum_{p\in Q_r}\mu(C_p)\le 1$. We find finite subsets $Q'_r\subseteq Q_r$ with $\sum_{p\in Q_r\setminus Q'_r}\mu(C_p)<\e_r$.

Consider the sets $C=\bigcup C_p$, where the union is taken for $p\in\bigcup_r(Q_r\setminus Q'_r)$ and let $Y=A\setminus C$. We have $\mu(C)<\sum_r\e_r$ and $\mu(Y)\ge\mu(A)-\mu(C)>1-\e$. As a closed and connected subset of $T$, $Y$ is subtree. It remains to prove that $Y$ is compact.

As a subtree $Y$ is a complete metric space, so it is enough to show that for any $\e^*>0$ and for any infinite sequence of points $y_i\in Y$ we can find $i\ne j$ with $d(y_i,y_j)<\e^*$. As the sequence $(d(p_0,y_i))_{i\in\N}$ is bounded we can find a short interval containing infinitely many of these distances, namely for some $a\ge0$ the subsequence consisting of the points $y_i$ satisfying $a\le d(p_0,y_i)\le a+\e^*/4$ is infinite. If $a=0$, then the distance of any two of the points in the subsequence is at most $\e^*/2$ and we are done. Otherwise, we find a positive rational number $r$ with $a-\e^*/4<r<a$. For a point $y_i$ in the subsequence the unique point $z_i$ in $[p_0,y_i]\cap Q_r$ satisfies $z_i\in Q'_r$ and $d(z_i,y_i)<\e^*/2$. As $Q'_r$ is finite we can find two distinct elements of the subsequence $y_i$ and $y_j$ with $z_i=z_j$. But then $d(y_i,y_j)\le d(y_i,z_i)+d(y_j,z_j)<\e^*$ as needed. \qed

The next lemma extends the concept of compact exhaustion to the domains of long dendrons. The similar statement for dendrons is even more immediate as we can always choose $a=1/2$.

\begin{lemma}\label{dencomp}
For a long dendron $D=(T,d,\nu)$ the real tree $(T,d)$ is separable, moreover it has a countable subset that is dense in every proper segment. For every $\e>0$ we have a compact subtree $Y$ of $T$ and $a>0$ such that $\nu(Y\times[0,a])>1-\e$.
\end{lemma} 

\proof
Recall that $\nu$ is a Borel probability measure on $A_D=T\times[0,\infty)$, and consider its marginal $\nu_T$ on $T$. In other words, $\nu_T$ is the push-forward of $\nu$ along the projection of $A_D$ to $T$. This $(T,d,\nu_T)$ a measured real tree with each branch having positive measure. We have a countable set in $T$ that is dense in every proper segment by Lemma~\ref{compact}. The lemma further implies that there is a compact subtree $Y$ of $T$ such that $\nu_T(Y)>1-\e/2$. We have $\lim_{a\to\infty}\nu(Y\times[0,a])=\nu(Y\times[0,\infty))=\nu_T(Y)>1-\e/2$. Therefore $a>0$ can be chosen as required. \qed

\begin{lemma}\label{minimal}
The unique minimal subtree of a real tree $T$ containing the points $x_1,\dots,x_n\in T$ is $\bigcup_{i=1}^n[x_1,x_i]$.
\end{lemma}

\proof
Any subset of $T$ containing $x_1$ and $x_i$ must also contain $[x_1,x_i]$ to be connected. But $\bigcup_{i=1}^n[x_1,x_i]$ is connected, closed and non-empty, so by Lemma~\ref{subtree}(1) it is indeed a subtree. \qed

\begin{definition}
Let $D=(T,d,\nu)$ be a long dendron and $r\ge1$ integer. We call $x\in A_D^n$ an {\bf$n$-sample} of $D$. For an $n$-sample $x=((p_1,a_1),\dots,(p_n,a_n))$, we define the measured real tree $T^D_x=(T',d',\mu)$ as follows. Let $T_0$ be the minimal subtree of $T$ containing the points $p_1,\dots,p_n$, whose existence is given by Lemma~\ref{minimal}. We obtain $T'$ from $T_0$ by appending the segments $[p_i,q_i]$ of length $a_i$ to $T_0$ for $i=1,\dots,n$ such that $T_0$ and all the sets $[p_i,q_i]\setminus\{p_i\}$ are pairwise disjoint. We set $d'$ to be the shortest path metric. This makes $(T',d')$ a real tree. The measure $\mu$ is defined by $\mu(H)=|\{1\le i\le n\mid q_i\in H\}|/n$ for any Borel subset $H\subseteq T'$. In other words, $\mu$ is the distribution of $q_i$ with a uniform random $i$. We also define the map $\alpha^D_x:T'\to A_D$ by setting $\alpha^D_x(p)=(p,0)$ for $p\in T_0$ and $\alpha^D_x(p)=(p_i,d'(p,p_i))$ if $p\in[p_i,q_i]$.

We call the sequence $(x^n)_{n\in\N}$ an {\bf infinite sample} of $D$ if $x^n$ is an $n$-sample for all $n$. An infinite sample obtained by independently selecting an $n$-sample $x^n$ for all $n$ according the distribution $\nu^n$ is called an {\bf infinite random sample} of $D$. We say that the infinite sample $(x^n)_{n\in\N}$ {\bf obeys} a Borel set $H\subseteq A_D$ if $\nu(H)=\lim_{n\to\infty}|\{1\le i\le n\mid x^n_i\in H\}|/n$, where $x^n=(x^n_1,\dots,x^n_n)$.
\end{definition}

\begin{lemma}\label{lln}
For a long dendron $D$ and every single Borel subset $H\subseteq A_D$ an infinite random sample of $D$ almost surely obeys $H$.
\end{lemma}

\proof
This is a form of the law of large numbers. \qed

Our main result in this section is the following theorem. Theorem~\ref{unique} will be a simple consequence.

\begin{theorem}\label{infsamp}
Let $(x^n)_{n\in\N}$ be an infinite random sample of the long dendron $D=(T,d,\nu)$. The measured real trees $T^D_{x^n}$ almost surely form an essentially bounded sequence. Furthermore, the long dendron associated with their ultraproduct is almost surely isomorphic to $D$.
\end{theorem}

The proof is through a series of lemmas. For this proof we fix $D=(T,d,\nu)$, and the infinite sample $(x^n)_{n\in\N}$ of $D$ with $x^n=(x^n_1,\dots,x^n_n)$. We say that a Borel set in $A_D$ is {\bf obeyed} if the sequence $(x^n)_{n\in\N}$ obeys it. Throughout the proof we will assume that various Borel sets in $A_D$ are obeyed. We can do that as long as we make the assumption for a countable family of sets by Lemma~\ref{lln}.

We introduce some notation. Let $T^D_{x^n}=(T_n,d_n,\mu_n)$ and let $\alpha_n=\alpha^D_{x^n}$. Let $(\bT,\bd)$ be the metric ultraproduct of $(T_n,d_n)$ and let $\hmu$ be the push-forward measure of $\prodo\mu_n$ along the natural projection $\hpi:\uT\to\bT$ from the set theoretic ultraproduct $\uT=\prodo T_n$ to the metric ultraproduct $\bT$.

Let $\alpha=\limo\alpha_n$. Here $\alpha_n:T_n\to A_D$ and $A_D$ (considered always with the product topology) is not compact, so $\alpha$ is defined on a subset of $\uT$, namely for $[(x_n)_{n\in\N}]\in\uT$, where $\limo\alpha_n(x_n)$ exists.

Consider the distance function $d'$ on $A_D$ defined as
\begin{equation}\label{d'}
d'((u,a),(v,b))=\left\{\begin{array}{ll}|a-b|&\hbox{ if }u=v\\d(u,v)+a+b&\hbox{ otherwise.}\end{array}\right.\end{equation}
Note that this is indeed a distance function making $(A_D,d')$ a real tree as mentioned in Remark~\ref{meareal}.

For any $n$ and $p,q\in T_n$ we have
$$d'(\alpha_n(p),\alpha_n(q))\le d_n(p,q)\le d_D(\alpha_n(p),\alpha_n(q))\,.$$
Note that $d_D$ is continuous on $A_D^2$ and $d'$ is lower semicontinuous there. Thus, for points $\up$ and $\uq$ in the domain of $\alpha$ we have
$$d'(\alpha(\up),\alpha(\uq))\le\ud(\up,\uq)\le d_D(\alpha(\up),\alpha(\uq))\,,$$
where $\ud=\limo d_n$. This means, in particular, that $\alpha(\up)=\alpha(\uq)$ if $\ud(\up,\uq)=0$. One can also see that if $\ud(\up,\uq)=0$ and one of $\alpha(\up)$ and $\alpha(\uq)$ is defined then so is the other. Therefore, we can define the function $\halpha$ by setting $\halpha(\hpi(\ux))=\alpha(\ux)$ if $\alpha(\ux)$ is defined and keeping $\halpha(\hpi(\ux))$ undefined if $\alpha(\ux)$ is not defined. Let $\bD$ stand for the domain of $\halpha$. For $\bp,\bq\in\bD$ we have
\begin{equation}\label{hbound}
d'(\halpha(\bp),\halpha(\bq))\le\bd(\bp,\bq)\le d_D(\halpha(\bp),\halpha(\bq))\,.
\end{equation}

We call a point $\bp\in\bD$ a {\bf base point} for $p\in T$ if $\halpha(\bp)=(p,0)$. The {\bf base} of $\bT$ is the set of base points in $\bD$.

Let us fix a countable set $Q\subseteq T$ that is dense in every proper segment of $T$. We can do this by Lemma~\ref{dencomp}. In the sequel we assume that $B\times[0,\infty)$ is obeyed if $B$ is a $q$-branch of $T$ for some $q\in Q$.  Note that $T$ is separable and $Q$ is countable, therefore we are making this assumption about a countable family of sets.

\begin{lemma}\label{base}
There is a unique base point for every $p\in T$. The function $\beta:T\to\bT$ mapping $p\in T$ to the base point for $p$ is an isometry from $T$ to the base of $\bT$.
\end{lemma}

\proof
Notice that $d'((p,0),(q,0))=d_D((p,0),(q,0))=d(p,q)$ for $p,q\in T$, so if $\bp$ is a base point for $p\in T$ and $\bq$ is a base point for $q\in T$, then
$$\bd(\bp,\bq)=d(p,q)$$
by our bound (\ref{hbound}). In particular, this means that the base point for $p\in T$ is unique if exists and if they exist for all $p\in T$, then $\beta$ is an isometry as claimed. It remains to show the existence.

Notice that $T\cap T_n$ is a subtree of $T$. Let $p$ be an arbitrary point of $T$ and let $\bp=\widehat{[(p_n)_{n\in\N}]}\in\bT$ with $p_n=\pi^{T,d}_{T\cap T_n}(p)$. Assume $\bp$ is not a base point for $p$. We have $\alpha_n(p_n)=(p_n,0)$, so this means that there is $\e>0$ such that $d(p,p_n)<\e$  for $\omega$-few indices $n$. In this case we can take a segment $[p,p']$ of length $\e$ in $T$, an intermediate point $q\in Q$ in that segment and notice that $q\in T_n$ for $\omega$-few indices $n$. As $q$ is an intermediate point, there are at least two distinct $q$-branches $B$ and $B'$. The sets $B\times[0,\infty)$ and $B'\times[0,\infty)$ are both obeyed and have positive $\nu$-measure. Therefore either is avoided by the sample $x^n$ for finitely many indices $n$. If neither is avoided, then $q\in T_n$. So $q\in T_n$ for $\omega$-few indices $n$, but $q\notin T_n$ for finitely many indices $n$. The contradiction shows that $\bp$ is a base point for $p$. \qed

Let us choose compact subtrees $Y_s$ of $T$ and reals $a_s>0$ for $s\in\N$ and form the compact product spaces $Z_s=Y_s\times[0,a_s]\subseteq A_D$. We can do this such that $\nu(Z_s)>1-1/s$ by Lemma~\ref{dencomp}. From now on we assume that the sets $Z_s$ are obeyed.

For $s,n\ge1$ let $T_{s,n}=\alpha_n^{-1}(Z_s)$. This is a subtree of $T_n$ whenever non-empty. Let $\bT_s=\hprodo T_{s,n}\subseteq\bT$.

\begin{lemma}\label{essential}
The sequence of metric measure spaces $(T^D_{x^n})_{n\in\N}$ is essentially bounded with $\bD$ contained in the dominant cluster of $\bT$.

The sets $\bT_s$ are real trees with $\hmu(\bT_s)\ge1-1/s$. We have $\bT_s\subseteq\bD$ and therefore $\hmu(\bD)=1$.
\end{lemma}

\proof
We have $\diam(T_{s,n})\le\diam(Y_s)+2a_s$. This makes the sequence of measured real trees $(T_{s,n},d_n)_{n\in\N}$ uniformly bounded and thus their metric ultraproduct $\bT_s$ is a real tree or empty by Lemma~\ref{treeprod}. We have $\hmu(\bT_s)\ge\limo\mu_n(T_{s,n})$. Here $\mu_n(T_{s,n})=|\{1\le i\le n\mid x^n_i\in Z_s\}/n$, so as $Z_s$ is obeyed, we have $\lim_{n\to\infty}\mu_n(T_{s,n})=\nu(Z_s)>1-1/s$. Therefore $\hmu(\bT_s)>1-1/s$. In particular, $\bT_s$ is not empty, so it is a real tree.

We have $\bT_s\subseteq\bD$ as for $\bp=\widehat{[(p_n)_{n\in\N}]}$ with $p_n\in T_{s,n}$ the ultralimit $\halpha(\bp)=\limo\alpha_n(p_n)$ exists because $\alpha_n(p_n)\in Z_s$ and $Z_s$ is compact. Thus, $\hmu(\bD)=1$ as claimed.

We have $\bd(\bp,\bq)<\infty$ for $\bp,\bq\in\bD$ by the bound (\ref{hbound}), therefore $\bD$ is contained in a single cluster of $\bT$. As $\hmu(\bD)=1$, this cluster must be dominant and the sequence $(T^D_{x^n})_{n\in\N}$ is essentially bounded as claimed. \qed

We denote the dominant cluster of $\bT$ by $\bX$. We slightly abuse notation by also writing $\bX$ when referring to the ultraproduct of the measured real trees $T^D_{x^n}$, that is, to $\bX$ together with the restrictions of $\bd$ and $\hmu$ making it a semi-measured real tree.

\begin{lemma}\label{corebase}
The core of $\bX$ is the base of $\bT$.
\end{lemma}

\proof
First we claim that every inner point $\bx$ of $\bX$ belongs to $\bD$. Indeed, for $\bx\in\bX\setminus\bD$ the real tree $\bT_s$ does not contain $\bx$, so it must be contained in a single $\bx$-branch. This branch has measure over $1-1/s$. As such an $\bx$-branch exists for each $s$, $\bx$ is not an inner point.

Consider a point $\bx=\widehat{[(x_n)_{n\in\N}]}\in\bD$ with $\halpha(\bx)=(z,a)$. Let $\alpha_n(x_n)=(z_n,a_n)$. We define a subtree $T'_n$ of $T_n$ for $n\in\N$ as follows. Recall that $T_n$ is obtained from a subtree of $T$ by appending $n$ segments. If $x_n\notin T$ we obtain $T'_n$ from $T_n$ by removing the added segment that contains $x_n$. If $x_n\in T$ (that is, $a_n=0$) we simply set $T'_n=T_n$. The distance between $x_n$ and $T'_n$ is $a_n$. Let $\bT'=\bX\cap\hprodo T'_n$. This is a subtree by Lemma~\ref{treeprod}. The distance between $\bx$ and $\bT'$ is $\limo a_n=a$, so if $a>0$, then $\bx\notin\bT'$ and $\bT'$ must lie inside a single $\bx$-branch. As $\mu_n(T'_n)\ge1-1/n$, we have $\hmu(\bT')=1$, so this $\bx$-branch has full measure and therefore $\bx$ is not an inner point if $a>0$.

So far we proved that for an inner point $\bx$ of $\bX$, $\halpha(\bx)$ must be defined and $\halpha(\bx)=(z,a)$ with $a>0$ is not possible. Therefore, every inner point of $\bX$ must be in the base.

Let $q\in Q$ be an intermediate point in $T$. Let $B$ and $B'$ be two distinct $q$-branches of $T$. Let $C=B\times[0,\infty)$ and $C'=B'\times[0,\infty)$. These positive measure subsets of $A_D$ are both obeyed. Therefore the sets $\bC=\hprodo C_n$ and $\bC'=\hprodo C'_n$ are also positive measure subsets, where $C_n=\alpha_n^{-1}(C)$ and $C'_n=\alpha_n^{-1}(C')$. Consider points $\bp=\widehat{[(p_n)_{n\in\N}]}\in\bC\cap\bX$ with $p_n\in C_n$ and $\bp'=\widehat{[(p'_n)_{n\in\N}]}\in\bC'\cap\bX$ with $p'_n\in C'_n$. We have $q\in[p_n,p'_n]$. Note that $p_n$ or $p'_n$ might be undefined for $\omega$-few indices $n$ where $C_n$ or $C'_n$ is empty, but otherwise $q\in T_n$, so we can set $\bq=\widehat{[(q_n)_{n\in\N}]}$ with $q_n=q$ if $q\in T_n$ and $q_n$ undefined otherwise. Clearly, $\bq$ is the base point for $q$. By Lemma~\ref{treeprod} we have $[\bp,\bp']=\hprodo[p_n,p'_n]$, so $\bq\in[\bp,\bp']$. This implies that $\bq$ is an inner point in $\bX$. Indeed, otherwise there would be a full measure $\bq$-branch $\bB$ and the positive measure sets $\bC$ and $\bC'$ would both intersect $\bB$, but the segment between any two points of $\bB$ does not contain $\bq$.

Assume $|T|>1$. Then the intermediate points in $Q$ are dense in $T$, so by the isometry, the base points for these, all inner points, are dense in the base. This means that $\core(\bX)$ contains the base. But we also saw that all inner points are in the base, and the base is complete, therefore closed, so we must have that $\core(\bX)$ is the base as claimed.

The same conclusion is even easier to obtain if $|T|=1$. Indeed, we know that the core is not empty (Lemma~\ref{core1}) and there is no inner point outside singleton base, so the core of $\bX$ must again be the base. \qed

Let $D^*=(T^*,d^*,\nu^*)$ be the long dendron associated with $\bX$ and let $\alpha^*:\bX\to A_{D^*}$ be the associated projection. Recall that $\beta$ is the function that maps a point in $T$ to the unique base point for it in $\bT$.

\begin{lemma}\label{alphastar}
Let $\bx\in\bD$ be a point with $\halpha(\bx)=(p,a)\in A_D$. Then $\alpha^*(\bx)=(\beta(p),a)\in A_{D^*}$.
\end{lemma}

\proof
Let $\bq\in\core(\bX)$, so by Lemma~\ref{corebase} $\bq$ is the base point for some $q\in T$. We have $$d'(\halpha(\bx),\halpha(\bq))\le\bd(\bx,\bq)\le d_D(\halpha(\bx),\halpha(\bq))$$
by bound (\ref{hbound}). Here $\halpha(\bx)=(p,a)$, $\halpha(\bq)=(q,0)$ and $d'((p,a),(q,0))=d_D((p,a),(q,0))=d(p,q)+a$. So we have $\bd(\bx,\bq)=d(p,q)+a$. This implies that $\pi^\bX(\bx)$, the unique closest point to $\bx$ in $\core(\bX)$ is $\beta(p)$ with $\bd(\bx,\beta(p))=a$. This proves the lemma. \qed

\begin{remark}
This last lemma gives a nice description of $\alpha^*$ on the domain $\bD$ of $\halpha$. It is somewhat harder to describe $\alpha^*$ on the zero measure set $\bX\setminus\bD$. For example for a point $\bp=\widehat{[(p_n)_{n\in\N}]}\in\bX\setminus\bD$ the second coordinate of $\alpha^*(\bp)$ is strictly larger than the ultralimit of the second coordinates of $\alpha_n(p_n)$. Fortunately, knowing $\alpha^*$ almost everywhere in $\bX$ is enough to determine $\nu^*$.
\end{remark}

\proof[Proof of Theorem~\ref{infsamp}.]
Let $(x^n)_{n\in\N}$ be an infinite random sample of the long dendron $D$. We will use the notation introduced above. By Lemma~\ref{lln} our assumption on various sets being obeyed holds almost surely, so we can use Lemmas~\ref{base}--\ref{alphastar} above.

In particular, by Lemma~\ref{essential} the sequence of measured real trees $(T^D_{x^n})_{n\in\N}$ is essentially bounded, so its ultraproduct $\bX$ exists. For the long dendron $D^*=(T^*,d^*,\nu^*)$ associated with $\bX$ we have $T^*=\core(\bX)$, so by Lemmas~\ref{base} and \ref{corebase} the function $\beta:T\to T^*$ is an isometry. It remains to prove that the function $\beta':A_D\to A_{D^*}$ defined by $\beta'(p,a)=(\beta(p),a)$ is almost surely measure preserving from $(A_D,\nu)$ to $(A_{D^*},\nu^*)$. That is, we need to show that almost surely
$$\nu(\beta'^{-1}(H))=\nu^*(H)$$
holds for all Borel sets $H\subseteq A_{D^*}$. Here $\nu^*$ is defined as the push-forward of $\hmu$ along $\alpha^*$. By Lemmas~\ref{alphastar} and \ref{essential} we have $\alpha^*(\bx)=\beta'(\halpha(\bx))$ almost everywhere in $\bX$. So we have 
$$\nu^*(H)=\hmu((\alpha^*)^{-1}(H))=\hmu(\halpha^{-1}(\beta'^{-1}(H)))\,.$$
Therefore, it is enough to show that almost surely
\begin{equation}\label{muequal}
\nu(H)=\hmu(\halpha^{-1}(H))
\end{equation}
holds for all Borel sets $H\subseteq A_D$.
Both sides of Equation~\ref{muequal} are Borel probability measures on $A_D$. It is therefore enough to check the equation for an intersection-closed family $\cA$ of Borel sets that generate the entire Borel $\sigma$-algebra on $A_D$. We choose such a countable family $\cA$ consisting of closed sets. For example if $T$ has more than a single point, then the finite intersections of sets of the form $(T\setminus B)\times[0,a]$ will work if $B$ runs through $q$-branches for $q\in Q$ and $a\ge 0$ is a rational number. Recall, that $Q$ is a countable set in the separable real tree $T$, so this is a countable family. Therefore, it is enough to check that Equation~\ref{muequal} holds almost surely for any single closed set $H\subseteq A_D$.

First we claim that the inequality
\begin{equation}\label{muless}
\nu(H)\le\hmu(\halpha^{-1}(H))
\end{equation}
holds for any closed set $H\subseteq A_D$ that $(x^n)_{n\in\N}$ obeys. Let us fix such a set $H$.
Consider the sets $H_n=\alpha_n^{-1}(H)$ for $n\in\N$. We have $\mu_n(H_n)=|\{1\le i\le n\mid x^n_i\in H\}|/n$ and $\lim_{n\to\infty}\mu_n(H_n)=\nu(H)$ as $H$ is obeyed. Consider $\bH=\hprodo H_n$. We have $\hmu(\bH)\ge\limo\mu_n(H_n)=\nu(H)$. Any point $\bp\in\bH\cap\bD$ can be written as $\bp=\widehat{[(p_n)_{n\in\N}]}$ with $\alpha_n(p_n)\in H$ and we have $\halpha(\bp)=\limo\alpha_n(p_n)$. An ultralimit of points in $H$ is also in $H$ as $H$ is closed. So we have $\bH\cap\bD\subseteq\halpha^{-1}(H)$. Using Lemma~\ref{essential}, we have $\nu(H)\le\hmu(\bH)=\hmu(\bH\cap\bD)\le\hmu(\halpha^{-1}(H))$ as claimed.

Let us now fix a closed set $H\subseteq A_D$. Our goal is to prove that Equation~\ref{muequal} holds almost surely. Let us choose an increasing sequence of closed sets $H'_1\subseteq H'_2\subseteq\cdots$ with $\bigcup_{s=1}^\infty H'_s=A_D\setminus H$. One can, for example, take $H'_s$ to be the complement of the open $1/s$-neighborhood of $H$ in some metrization of $A_D$. We have $\lim_{s\to\infty}\nu(H'_s)=\nu(A_D\setminus H)=1-\nu(H)$ and similarly, using that $\halpha$ is defined almost everywhere (Lemma~\ref{essential}), we have $\lim_{s\to\infty}\hmu(\halpha^{-1}(H'_s))=\hmu(\halpha^{-1}(A_D\setminus H))=1-\hmu(\halpha^{-1}(H))$.

By lemma~\ref{lln} $H$ and all the sets $H'_s$ are almost surely obeyed. If so, we have $\nu(H)\le\hmu(\halpha^{-1}(H))$ and also $\nu(H'_s)\le\hmu(\halpha^{-1}(H'_s))$ for all $s$, all from Inequality~(\ref{muless}). Taking limit on both sides of the latter inequality we get $1-\nu(H)\le1-\hmu(\halpha^{-1}(H))$. This makes Equation~(\ref{muequal}) hold almost surely and finishes the proof of the theorem. \qed

\begin{definition}
Let $n\ge1$ and let $A\in M_n$ a $n$ by $n$ real matrix. The measured real tree $(T,d,\mu)$ is an {\bf$A$-tree} if there exist $q=(q_1,\dots,q_n)\in T^n$ such that (i) $\rho^{T,d}_n(q)=A$, (ii) no subtree of $T$ other than $T$ itself contains all of the points $q_i$ and (iii) $\mu$ is a Borel measure on $T$ given by $\mu(H)=|\{1\le i\le n\mid q_i\in H\}|/n$.
\end{definition}

\begin{lemma}\label{atree}
Let $n\ge1$ and $A\in M_n$. If there exists an $A$-tree it is unique up to measure preserving isometry. If $n\ge2$ and $x$ is an $n$-sample of a long dendron $D$, then $T^D_x$ is a $\rho^{A_D,d_D}_n(x)$-tree.
\end{lemma}

\proof
Assume both $(T,d,\mu)$ and $(T',d',\mu')$ are $A$-trees for $A=(d_{i,j})_{1\le i,j\le n}$. Let $q=(q_1,\dots,q_n)\in T^n$ and $q'=(q'_1,\dots,q'_n)\in T'^n$ be the samples that satisfy conditions (i-iii) in the definition. By condition (ii) and Lemma~\ref{minimal} we have $T=\bigcup_{i=1}^n[q_1,q_i]$ and $T'=\bigcup_{i=1}^n[q'_1,q'_n]$. Both the segments $[q_1,q_i]$ and $[q'_1,q'_i]$ have length $d_{1,i}$ by condition (i), so there is an isometry $f_i$ between them with $f_i(q_1)=q'_1$, $f_i(q_i)=q'_i$. The intersection of the segments $[q_1,q_i]$ and $[q_1,q_j]$ is a segment starting at $q_1$, and by condition (i) its length is $(d_{1,i}+d_{1,j}-d_{i,j})/2$. Similarly, $[q'_1,q'_i]\cap[q'_1,q'_j]$ is a segment of the same length starting at $q'_1$. Therefore, the functions $f_i$ agree on the intersections of their domains, so there is a global function $f:T\to T'$ extending all these functions. It is easy to see that $f$ is an isometry between $T$ and $T'$. By condition (iii) and since $f(q_i)=q'_i$ for all $i$, the map $f$ is measure preserving.

Let us recall, that for a long dendron $D=(T,d,\nu)$ and a $n$-sample $x=((p_1,a_1),\dots,(p_n,a_n))$ of $D$ we constructed the measured real tree $T^D_x$  by adding segments $[p_i,q_i]$ to the minimal subtree of $T$ containing the points $p_i$. The measured real tree $T^D_x$ with $q=(q_1,\dots,q_n)$ satisfy conditions (i) and (iii) for a $\rho^{A_D,d_D}_n(x)$-tree. They also satisfy condition (ii) as we assumed $n\ge2$. Note that condition (ii) might be violated for $n=1$ as in that case $T^D_x$ is a segment of length $a_1$ and (if $a_1>0$) not the trivial subtree $\{q_1\}$. \qed. 

\proof[Proof of Theorem~\ref{unique}.]
Using the observations in Lemma~\ref{atree} we can restate Theorem~\ref{infsamp} as follows. For a long dendron $D$ take independent samples $A_n$ from the distribution $\tau_n(D)$ and create the corresponding $A_n$-trees $T_n$. These exist and are unique up to measure preserving isometry. Almost surely, the sequence $(T_n)_{n\in\N}$ is essentially bounded and the long dendron associated with its ultraproduct is isomorphic to $D$.

If the two long dendrons $D$ and $D'$ satisfy $\tau_n(D)=\tau_n(D')$ for all $n$, then the above process is the same for both of them. Almost surely, the process results in an associated long dendron that is isomorphic to both $D$ and $D'$. Thus, $D$ and $D'$ are isomorphic as claimed.
\qed

\section{The proof of Theorem \ref{tetel2}}\label{pr2}

\begin{definition} Let $T$ be a finite (graph-theoretic) tree. We turn it into a real tree by keeping the vertices and replacing every edge by an arbitrary length interval connecting the corresponding vertices. Finally, we make it into a measured real tree by adding an arbitrary  Borel measure concentrated on the original set of vertices. We call the measured real trees obtained this way {\bf finite real trees}. Observe that a measured real tree $(T,d,\mu)$ is a finite real tree if and only if (i) there is a finite set of points in $T$ with no subtree containing all of them other than $T$ and (ii) $\mu$ is concentrated on a finite set of points. In particular,  if $x$ is an $n$-sample of a long dendron $D$, then $T^D_x$ is a finite real tree.
\end{definition}

Theorem~\ref{tetel2} easily follows from the following two lemmas as the weak
topology of $\prob(M_r)$ is metrizable for every $r\in\N$, by
Prokhorov's Theorem. Note that Lemma~\ref{samplelimit} has a direct proof using the Azuma Inequality (see 5.3 \cite{LSZ}), but at this point it is easier for us to use Theorems~\ref{tetel1+} and \ref{infsamp} instead.

\begin{lemma}\label{samplelimit}
Any long dendron $D$ is the sampling limit of finite real trees, that is, one can find finite real trees $(T_n,d_n,\mu_n)$ such that for every $r\ge1$ the sampling measures $\tau_r(T_n,d_n,\mu_n)$ weakly converge to $\tau_r(D)$.

If $D$ is a dendron, then the finite real trees can be chosen such that the diameter of each is at most one.
\end{lemma}

\proof
Let us choose an infinite random sample $(x^n)_{n\in\N}$ of $D$. By Theorem~\ref{infsamp}, $D$ can almost surely be obtained as the long dendron associated with the ultraproduct of the essentially bounded sequence of finite real trees $T^D_{x^n}$. In this case we have $\limo\tau_r(T^D_{x_n})=\tau_r(D)$ for any $r\in\N$ by Theorem~\ref{tetel1+}. The ultralimit is always the (ordinary) limit of a subsequence in a metrizable space, so we have a sequence of finite real trees $(T_n,d_n,\mu_n)$ satisfying that $\tau_r(T_n,d_n,\mu_n)$ weakly converges to $\tau_r(D)$. We can clearly find the same sequence  for all $r$.

If $D$ is a dendron, then $\diam(T^D_{x^n})\le1$ holds almost surely, so the second statement of the lemma also holds. \qed

\begin{lemma}
For a finite real tree $(T,d,\mu)$ of diameter at most $1$ one has a sequence of finite (graph-theoretic) trees $T_n$ such that $\tau_r(T_n)$ weakly converges to $\tau_r(T,d,\mu)$ for every $r\in\N$.
\end{lemma}

\proof
Let $(T,d)$ be obtained from the finite tree $(V,E)$ by replacing each edge $e\in E$ with a segment of length $a_e$. Note here that we have $a_e\le1$. For $n\in\N$ we build $T_n$ in three steps. In the first step we replace each edge $e$ of the tree $(V,E)$ with a path of length $\lceil a_en\rceil$. Then for every vertex $v\in V$ we form a set $H^v_n$ of $\lceil\mu(v)n^2\rceil$ new vertices and add them to the tree as leaves whose only neighbor is $v$. Finally, if the diameter of the tree constructed so far is below $n$, we attach a path to it that brings the diameter to exactly $n$.

The diameter of the tree $T_n$ so constructed is $n+O(1)$. The distance between any vertex $x\in H^v_n$ and $y\in H^w_n$ , is $d(v,w)n+O(1)$, so the normalized distance $d_n(x,y)$ is $d(v,w)+O(1/n)$. $T_n$ has $n^2+O(n)$ vertices, so the probability that a uniformly random vertex falls in $H^v_n$ is $\mu(v)+O(1/n)$ with the probability of choosing a vertex outside all the sets $H^v_n$ being $O(1/n)$. The hidden constants in the order notation depends only on the size of the vertex set $V$.

This shows that a $\mu$-random point $v$ can be coupled with a uniform random vertex of $T$ such that a point $v\in V$ is coupled with a vertex $x\in H^v_n$ with probability $1-O(1/n)$. Doing this in every coordinate one can couple a sample $(v_1,\dots,v_r)$ of $\mu^r$ with a uniform sample $(x_1,\dots,x_r)$ of $(V(T_n))^r$ such that $x_i\in H^{v_i}_n$ for all $1\le i\le r$ with probability $1-O(r/n)$. If this happens, then the matrices $\rho_r^{T,d}(v_1,\dots,v_r)$ and $\rho^{T_n,d_n}_r(x_1,\dots,x_n)$ differ by $O(1/n)$ in every coordinate. Therefore, their distributions $\tau_r(T,d,\mu)$ and $\tau_r(T_n)$ are close. In particular, this implies the weak convergence stated in the lemma. \qed

\end{document}